# ATTRACTION TIME FOR STRONGLY REINFORCED WALKS

By Codina Cotar and Vlada Limic[1]

*TU Berlin and CNRS*

We consider a class of strongly edge-reinforced random walks, where the corresponding reinforcement weight function is nondecreasing. It is known, from Limic and Tarrès [*Ann. Probab.* (2007), to appear], that the attracting edge emerges with probability 1 whenever the underlying graph is locally bounded. We study the asymptotic behavior of the tail distribution of the (random) time of attraction. In particular, we obtain exact (up to a multiplicative constant) asymptotics if the underlying graph has two edges. Next, we show some extensions in the setting of finite graphs, and infinite graphs with bounded degree. As a corollary, we obtain the fact that if the reinforcement weight has the form $w(k) = k^\rho$, $\rho > 1$, then (universally over finite graphs) the expected time to attraction is infinite if and only if $\rho \leq 1 + \frac{1+\sqrt{5}}{2}$.

**1. Introduction.** Let $\mathcal{G}$ be a locally finite graph with the edge set $E(\mathcal{G})$ and the vertex set $V(\mathcal{G})$. We will assume without further mention that $\mathcal{G}$ is connected. We call any two vertices $u, v$ connected by an edge *adjacent* (or *neighboring*); in this case, we write $u \sim v$ and denote by $\{u, v\} = \{v, u\}$ the edge connecting them. We will denote by

$$|\mathcal{G}| = |E(\mathcal{G})|$$

the number of edges of $\mathcal{G}$ and by

$$\#\mathcal{G} = |V(\mathcal{G})|$$

the number of vertices of $\mathcal{G}$. Finally, we denote by $D(\mathcal{G}) = \sup_{v \in V(\mathcal{G})} \text{degree}(v)$ the degree of $\mathcal{G}$, where, for any $v \in V(\mathcal{G})$, $\text{degree}(v)$ equals the number of edges incident to $v$.

Received December 2006; revised September 2008.

[1]Supported in part by an NSERC research grant and by the Alfred P. Sloan Research Fellowship.

Let $(\ell_0^e, e \in E(\mathcal{G}))$ be given integers and assume that $\ell_0^e \geq 0$, $e \in E(\mathcal{G})$. Given a *reinforcement weight function* $w \colon \{0, 1, 2, \ldots\} \mapsto (0, \infty)$, the edge-reinforced random walk (ERRW) on $\mathcal{G}$ records nearest neighbor step transitions of a particle in $V(\mathcal{G})$. That is:

(i) if currently at vertex $v \in V(\mathcal{G})$, in the next step, the particle jumps to a vertex $u \in V(\mathcal{G})$ adjacent to $v$;

(ii) the probability of a jump to $u$ is $w$-*proportional* to the number of previous traversals of the edge $\{v, u\}$.

The more formal definition is as follows. If $\mathcal{G}$ is a finite graph, it seems natural, from the point of notation, to construct and study the edge-reinforced random walk started at the *initial time*

$$t_0 := \sum_{e \in E(\mathcal{G})} \ell_0^e \geq 0;$$

a process starting at time $0$ is obtained by a time shift. If $\mathcal{G}$ is an infinite graph, we simply set $t_0 := 0$. Denote by $I_n$ the (random) position of the edge-reinforced random walk at time $n$. Then, $I_{t_0} \in V(\mathcal{G})$ is the initial position and $\{I_n, I_{n+1}\} \in E(\mathcal{G})$ for all $n \geq t_0$, almost surely. Let $\mathcal{F}_n$ be the filtration

(1) $$\mathcal{F}_n = \sigma\{I_k, k = 0, \ldots, n, (\ell_0^e, e \in E(\mathcal{G}))\}.$$

Moreover, the dynamics of the edge-reinforced random walk is prescribed according to the rule

$$P(I_{n+1} = v | \mathcal{F}_n) 1_{\{I_n = u\}} = \frac{w(X_n^{\{u,v\}})}{\sum_{y \sim u} w(X_n^{\{u,y\}})} 1_{\{I_n = u, u \sim v\}},$$

where, for any $e \in E(\mathcal{G})$,

(2) $$X_n^e = \ell_0^e + \sum_{i=t_0}^{n-1} 1_{\{e \text{ was traversed at } i\text{th step}\}} = \ell_0^e + \sum_{i=t_0}^{n-1} 1_{\{\{I_i, I_{i+1}\} = e\}}$$

equals the *initial weight* $\ell_0^e$ incremented by the total number of (undirected) traversals of edge $e$ prior to time $n$. Note that $t_0$ is chosen so that whenever $V(\mathcal{G}) < \infty$, $\sum_{e \in E(\mathcal{G})} X_k^e = k$ for all $k \geq t_0$, almost surely. The starting weights $X_{t_0} := \ell_0^e$ are specified as deterministic above, but one could use random variables instead in applications and definition (1) accounts for this possibility. Our results would then hold conditionally on the starting weights.

We denote by $\mathcal{G}_1$ the *range* of the edge-reinforced random walk on $\mathcal{G}$. More precisely, we let

$$\mathcal{G}_1 = (V(\mathcal{G}_1), E(\mathcal{G}_1))$$

be the random subgraph of $\mathcal{G}$ where, for any $v \in V(\mathcal{G})$, we have

$$v \in V(\mathcal{G}_1) \quad \Leftrightarrow \quad \exists n \geq t_0 \qquad \text{such that} \quad I_n = v,$$

and, for any $e \in E(\mathcal{G})$,

$$e \in E(\mathcal{G}_1) \quad \Leftrightarrow \quad \exists n \geq t_0 \qquad \text{such that} \quad \{I_n, I_{n+1}\} = e.$$

Apart from the behavior analogous to recurrence or transience of Markov chains (see, e.g., [7, 9] or [10], Theorems 5.2 and 5.6), ERRW may exhibit a very different asymptotic behavior as time increases. For example, it is easy to see, [5, 11], that the assumption

(A0) $$\sum_k \frac{1}{w(k)} < \infty$$

is sufficient for the event

$$\{\mathcal{G}_1 \text{ is a finite graph}\}$$

to have probability 1, whenever $D(\mathcal{G}) < \infty$. It is easy to find examples of locally bounded trees with $D(\mathcal{G}) = \infty$ such that (A0) holds but $\mathcal{G}_1$ is infinite with positive probability. Sellke [11] provides (slightly peculiar) examples of edge-reinforced random walks on $\mathbb{Z}$ where $\sum_k 1/w(k)$ is finite over even $k$ and infinite over odd $k$, but where $\mathcal{G}_1$ is still a finite graph, almost surely.

Next, we briefly discuss links between our work and the recent literature. For a detailed review of a number of interesting results on edge-reinforced random walks, we refer the reader to a recent survey of Pemantle [10] on stochastic reinforcement processes.

A result of Sellke [11] (the argument is also described in detail in [5], Section 2) implies that (A0) is sufficient and necessary for

(3) $$P(\text{the walk ultimately traverses a single edge}) = 1,$$

whenever the underlying graph is bipartite and of bounded degree. Limic [5] proves that (A0) implies (3) on any graph of bounded degree, where the reinforcement weight is a reciprocally summable power function. In a recent work, Limic and Tarrès [6] show that for a fairly general class of reinforcement weights [in particular, whenever $w$ is a nondecreasing function satisfying (A0)] (3) holds on any graph of bounded degree. We will refer to any weight $w$ satisfying condition (A0) as *strong* and to the corresponding ERRW as a *strongly reinforced walk*.

The current paper assumes the setting of [6] and is devoted to the study of the tail behavior of the *time of attraction*

(4) $$T = \inf\left\{k \geq 0 : \exists e \in E(\mathcal{G}) \text{ such that } \forall f \neq e\, X_k^f = \max_{m \geq k} X_m^f\right\}$$
$$= \inf\{k \geq 0 : \{I_n, I_{n+1}\} = \{I_{n+1}, I_{n+2}\}, \forall n \geq k\},$$

that is, the first time after which only the *attracting* edge is traversed. This random variable is an important statistic, useful for applications (e.g., [3] or [4]). In this paper, we make a few connections to the literature on the behavioral science of social insects and refer the reader to [10] for a diverse list of potential applications.

It is clear that the sequence of tail probabilities $(P(T > k), k \geq 1)$ depends on the structure of the underlying graph $\mathcal{G}$, the weight function $w$, the initial weights $\ell_0$ and the initial position $I_{t_0}$. However, the results of Sections 3.1–3.2 verify an interesting universality-type behavior. Namely, fix $w$ satisfying (A0), let $\mathcal{G}$ be an arbitrary finite graph with some prescribed initial edge weights and initial position, and let $\mathcal{G}'$ be the simple two-edge graph from Section 2 with initial weights equal to 1 on both edges. Then, if $P^{\mathcal{G}}$ (resp., $P^{1,1}$) denotes the law of the corresponding ERRW on $\mathcal{G}$ (resp., $\mathcal{G}'$), the asymptotic order of magnitude of $P^{\mathcal{G}}(T > k)$ is induced by that of $P^{1,1}(T > k)$. To some extent, this also holds on infinite trees of bounded degree; see Corollary 18.

DEFINITION 1. For sequences $a_k, b_k$ of real numbers, we write $a_k \asymp b_k$ if and only if

$$\frac{a_k}{b_k} \in [c, C], \qquad k \geq 1 \text{ for some } c, C \in (0, \infty).$$

The rest of the paper is organized as follows. Section 2 is devoted to the careful study of the two-edge setting. In particular, in Section 2.1, we prove Corollary 5, stated less precisely as follows. If we assume that $E(\mathcal{G})$ consists of two elements and let

$$Z_\infty := \#\{\text{times the nonattracting edge is traversed}\},$$

then, under assumption (A0),

$$P(Z_\infty = \ell) \asymp \frac{1}{w(\ell)}, \qquad \ell \to \infty.$$

A stronger statement, due to R. Pemantle (personal communication),

$$\text{(5)} \qquad \lim_{\ell \to \infty} w(\ell) P(Z_\infty = \ell) \in (0, \infty),$$

holds in this simple setting. The limit will become apparent in the course of the proof sketched at the end of Section 2.2.

Lemma 6 of Section 2.3 provides the initial order-of-*magnitude* estimates on the tail probability $P(T > k)$ as $k \to \infty$. More precisely, we prove that

$$P(T > n) \asymp \sum_{k=n+1}^{\infty} \sum_{\ell \leq k/2} \left[ \frac{1}{w(k-\ell)} \prod_{j=0}^{\infty} \frac{w(k+j-\ell)}{w(k+j-\ell) + w(\ell+1)} \right.$$
$$\left. + \frac{1}{w(\ell)} \prod_{j=0}^{\infty} \frac{w(j+\ell-1)}{w(j+\ell-1) + w(k+1-\ell)} \right].$$

Such an expression seems awkward for applications and we work further to find simplifications. In particular, Theorem 9 shows simpler looking asymptotics of the tail distribution of $T$, under the additional assumption (A1) that $w$ is a nondecreasing function. The main idea is simple: the event which with overwhelming probability contributes to the event $\{T = k + 1\}$ of interest is the one where at time $k$, the weaker edge (i.e., the edge with the lower current number of traversals) is traversed and at all future times, the stronger edge is traversed. Therefore ($Z_k$ denoting the number of traversals of the less traversed edge at time $k$),

$$P(T = k + 1) \asymp \sum_{\ell \leq k/2} P(Z_k = \ell) \frac{w(\ell)}{w(k - \ell) + w(\ell)} P(Z_\infty = \ell + 1 | Z_{k+1} = \ell + 1).$$

For $\ell$ close to $k/2$, it is plausible that $P(Z_k = \ell)$ is sufficiently small so that the contribution in the above sum vanishes asymptotically. For $\ell$ small, the middle term $w(\ell)/(w(k - \ell) + w(\ell))$ is again small. In order to estimate the above sum well, one then needs to find the interval of indices $\ell$ which make up the bulk of the contribution. As shown in Proposition 2, $P(Z_k = \ell) \asymp \frac{w(k-\ell)+w(\ell)}{w(k-\ell)w(\ell)}$, so it is plausible that the overwhelming contribution to the sum comes, approximately, from the range of indices where $P(Z_\infty = \ell + 1 | Z_{k+1} = \ell + 1) \asymp 1$. For formal estimates, see Section 2.3.

In Section 2.4, we include specific calculations for cases of $w$ that have already been used (or might be used) in applications (see [2, 4]) that satisfy the assumptions (A0)–(A1). In particular, we paraphrase as follows.

THEOREM 10(a).   *If $w(k) = k^\rho$ for some fixed $\rho > 1$ and if $\rho' = (\rho - 1)/\rho$, then*

$$P(T > k) \asymp \frac{1}{k^{\rho - \rho' - 1}}.$$

*In particular, $E(T)$ is infinite if $\rho \leq 1 + \frac{1+\sqrt{5}}{2}$ and finite if $\rho > 1 + \frac{1+\sqrt{5}}{2}$.*

This type of result should be particularly interesting for applications. In fact, in [4], for a similar model, the reinforcement weight is set to $w(k) = k^\rho$ and real-life data is compared to different values of $\rho$ and initial configurations. More precisely, the authors study a colony of ants which randomly explores a chemically unmarked territory, starting from its nest. The exploration is carried out on two branches $A$ and $B$. Initially, both branches are equally likely to be chosen. However, each ant that passes along one of the two branches leaves an additional pheromone mark and in this way influences the following ant's decision in choosing $A$ or $B$. In the real-life experiment, it is observed that, after initial fluctuations, one of the two branches becomes more or less completely preferred to the other. In their (reinforcement) model, $k$ represents the number of ants that have chosen a particular branch and $\rho$ determines the degree of nonlinearity. The model is used for further study of the explorer movement pattern in two-dimensional space.

Section 3 is devoted to analysis on general graphs of bounded degree. In particular, we are interested in a universality-type behavior of the tail distribution $P^{\mathcal{G}}(T > \cdot)$ over graphs once the reinforcement weight function $w$ is fixed. In the course of our analysis, we also obtain exponential bounds (see Lemma 25) on the tail distribution of $|\mathcal{G}_1|$ and, in particular, some information on the distance of the attracting edge from the starting point. Providing a universal lower bound on $P^{\mathcal{G}}(T > \cdot)$ in terms of the corresponding quantity in the two-edge graph setting turns out to be simple (see Lemma 12); however, finding an analogous upper bound is not as simple. Section 3.1 is devoted to analysis on trees. Here, initial universality-type behavior is demonstrated using comparison (coupling) arguments. Section 3.2 is devoted to the finite graph setting. By generalizing the technique of Section 2, a fairly general universality-type behavior is shown, under the additional assumption (A2). Finally, Section 3.3 discusses extensions to the infinite graph setting.

In the remainder of the paper, we assume that all edges have "*trivial*" *initial weight* $\ell_0 \equiv 1$, unless otherwise specified. Also, we will denote by $a \wedge b$ (resp., $a \vee b$) the minimum (resp., maximum) of two numbers $a$ and $b$, and by $\lfloor a \rfloor$, the integer part of a number $a$.

## 2. Two-edge case.

The ERRW on graph $\mathcal{G}$ that contains only two edges is the prototype model of interest. Several interesting qualitative features, specific to edge-reinforcement with particular reinforcement weight function $w$, are already observed and are usually relatively easy to verify. This process also corresponds to a generalized urn model; see, for example, [1] or [11]. A recent study by Oliveira and Spencer [8] concerns finer properties of this urn model in the case where $w(k) = k^\rho$ for some $\rho > 1$.

We will initially assume that $\mathcal{G}$ contains two vertices, 0 and 1, and two edges, *green* and *red*, connecting them. We abbreviate

$$G_n := X_n^{green}, \qquad R_n := X_n^{red}.$$

In the remainder of this section, we also assume that the initial configuration on the two edges is $G_2 = R_2 = 1$, unless otherwise specified. We use the notation $P^{a,b}$ for the law of the system with the initial configuration $G_{a+b} = a, R_{a+b} = b$. When there is no risk of confusion, we simply use $P$ for the law $P^{1,1}$.

The other natural choice of a graph with two edges is the one spanned by three vertices, $-1, 0$ and 1, with a green edge that connects 0 and $-1$ and a red one that connects 0 and 1. In the study of this model, we mainly concentrate on the case where the initial weights $a$ and $b$ are of opposite parity. We denote by $\bar{P}^{a,b}$ the law of the ERRW on the two-edge graph spanned by $-1, 0$ and $1$, started (without loss of generality) at the initial position $0$. Note that the study of $\bar{P}^{a,b}$ is necessary as it will be needed later for the subsequent analyses of the time of attraction of the ERRW on trees and on finite graphs. The main results in Section 3 are expressed in terms of both $P^{a,b}$ and $\bar{P}^{a,b}$.

Observe that under $\bar{P}^{a,b}$, we have $G_{a+b} = a, R_{a+b} = b$ and $G_{a+b+2j} - G_{a+b+2j-2} \in \{0, 2\}$, $R_{a+b+2j} - R_{a+b+2j-2} + G_{a+b+2j} - G_{a+b+2j-2} = 2$ for all $j \geq 1$.

2.1. *Some preliminary estimates.* Due to monotonicity, $R_\infty := \lim_{k\to\infty} R_k$ and $G_\infty := \lim_{k\to\infty} G_k$ exist almost surely as $(0, \infty]$-valued random variables. Define

$$Z_k = \min\{G_k, R_k\}.$$

Note that $Z_k \leq Z_{k+1}$ and that the limit

$$(6) \qquad Z_\infty := \lim_{k\to\infty} Z_k = R_\infty \wedge G_\infty$$

is an almost surely finite random variable since the reinforcement is strong.

PROPOSITION 2. *Define $c := \max_{k\geq 2}\left(\frac{w(k-1)+w(1)}{w(k-1)}\right)$.*

(a) *For any $1 \leq \ell < \frac{k}{2}$, we have*

$$\frac{w(1)}{c} P(Z_\infty = 1)\frac{w(k-\ell)+w(\ell)}{w(k-\ell)w(\ell)} \leq P(Z_k = \ell) \leq w(1)\frac{w(k-\ell)+w(\ell)}{w(k-\ell)w(\ell)}.$$

(b) *For any $\ell \geq 1$, we have*

$$\frac{w(1)}{c} \cdot \frac{P(Z_\infty = 1)}{w(\ell)} \leq P(Z_{2\ell} = \ell) \leq \frac{w(1)}{w(\ell)}.$$

Note that the lower bounds are interesting only for strongly reinforced walks, where $P(Z_\infty < \infty) = 1$ and $P(Z_\infty = 1) > 0$. A careful reader of the proof will note that all of the above inequalities are strict; however, we do not anticipate any use of this fact.

PROOF OF PROPOSITION 2. We will prove the upper bounds by induction and the lower bounds will follow in a similar way, as indicated at the end of the proof. First, note that for $\ell \leq \frac{k}{2} - 1$,

$$(7) \qquad \begin{aligned} P(Z_k = \ell) &= P(Z_{k-1} = \ell)\frac{w(k-\ell-1)}{w(k-\ell-1)+w(\ell)} \\ &\quad + P(Z_{k-1} = \ell-1)\frac{w(\ell-1)}{w(\ell-1)+w(k-\ell)}. \end{aligned}$$

Similarly, we also have, in the special case $k = 2\ell$,

$$(8) \qquad P(Z_{2\ell} = \ell) = P(Z_{2\ell-1} = \ell - 1)\frac{w(\ell - 1)}{w(\ell - 1) + w(\ell)}$$

and in the special case $k = 2\ell + 1$,

$$(9) \quad P(Z_{2\ell+1} = \ell) = P(Z_{2\ell} = \ell) + P(Z_{2\ell} = \ell - 1)\frac{w(\ell - 1)}{w(\ell - 1) + w(\ell + 1)}.$$

Since any probability is bounded by 1, we have, trivially,

$$P(Z_k = 1) < w(1)\frac{w(k - 1) + w(1)}{w(k - 1)w(1)},$$

an observation that will be used in the base and in each step of the induction.

The base of induction is the case $\ell = 1$, $k = 2\ell + 1 = 3$ and the statement here is trivial, as noted above.

Let us now assume that the upper bound inequalities in the statements (a) and (b) of the theorem hold for all $i \leq k - 1$ and $\ell \leq \frac{i}{2}$. For the induction step, we need to show that the bounds hold for $i = k$ and $\ell \leq \frac{k}{2}$.

Suppose, first, that $\ell \leq \lfloor \frac{k}{2} \rfloor - 1$. Then, by (7) and the induction hypothesis, we have

$$P(Z_k = \ell) \leq \frac{w(1)}{w(\ell)} + \frac{w(1)}{w(k - \ell)} = w(1)\frac{w(k - \ell) + w(\ell)}{w(k - \ell)w(\ell)}.$$

For the two atypical cases $k = 2\ell$ and $k = 2\ell + 1$, we have, similarly, by (8),

$$P(Z_{2\ell} = \ell) = P(Z_{2\ell-1} = \ell - 1)\frac{w(\ell - 1)}{w(\ell - 1) + w(\ell)} < \frac{w(1)}{w(\ell)}$$

and by (9),

$$P(Z_{2\ell+1} = \ell) < \frac{w(1)}{w(\ell)} + \frac{w(1)}{w(\ell + 1)} = w(1)\frac{w(\ell) + w(\ell + 1)}{w(\ell)w(\ell + 1)}.$$

The proof of the lower bounds is symmetric. Note that $P(Z_k = 1) \geq P(Z_\infty = 1)$ and the choice of $c$ was precisely made so that the lower bound holds both in (a) for any $k \geq 3$ and $\ell = 1$, and in (b) for $\ell = 1$. Given these initial bounds, the above argument by induction on $k$ will carry over to yield the lower bound of (a) and (b).  □

The result above under the law $P^{1,1}$ generalizes to the setting of the law $\bar{P}^{1,2}$ on a two-edge graph with three vertices, in the following way.

Proposition 3. *Define*

$$\bar{c} := \max_{k \geq 2 \text{ even}}\left(\frac{w(k) + w(1)}{w(k)}\right) \vee \max_{k \geq 1 \text{ odd}}\left(\frac{w(k) + w(2)}{w(k)}\right).$$

*For any $k \geq 1$, $1 \leq o, e \leq 2k$ such that $o$ is odd, $e$ is even and $o + e = 2k + 1$, we have*

$$\bar{P}^{1,2}(G_{2k+1} = o, R_{2k+1} = e) \leq (w(1) \vee w(2)) \frac{w(o) + w(e)}{w(o)w(e)}$$

*and*

$$\frac{w(1) \wedge w(2)}{\bar{c}} \cdot \frac{w(o) + w(e)}{w(o)w(e)} \prod_{j=1}^{2} \bar{P}^{1,2}(Z_\infty = j) \leq \bar{P}^{1,2}(G_{2k+1} = o, R_{2k+1} = e).$$

PROOF. We abbreviate $\bar{P}^{1,2}$ as $P$. First, we concentrate on the upper bound. If either $o = 1$ or $e = 2$ (or both), the upper bound is trivial, so the base of induction is verified. Now, as in the previous proposition, if both $o > 1$ and $e > 2$, we apply the induction step using

$$P(G_{2k} = o, R_{2k} = e) = P(G_{2k-2} = o - 2, R_{2k-2} = e) \frac{w(o-2)}{w(o-2) + w(e)}$$

$$+ P(G_{2k-2} = o, R_{2k-2} = e - 2) \frac{w(e-2)}{w(o-2) + w(e-2)}.$$

Similarly, note that $\bar{P}^{1,2}(G_{2k+1} = 1, R_{2k+1} = 2k) \geq \bar{P}^{1,2}(G_\infty = j) = \bar{P}^{1,2}(Z_\infty = j)$ and $\bar{P}^{1,2}(G_{2k+1} = 2k-1, R_{2k+1} = 2) \geq \bar{P}^{1,2}(R_\infty = 2) = \bar{P}^{1,2}(Z_\infty = 2)$, so the lower bound holds for any $k \geq 1$ whenever $o = 1$ or $e = 2$, with the above choice of $\bar{c}$. Given these initial bounds, one applies the induction step once again to prove the general lower bound. □

Moreover, using the same technique as above, one arrives at the following general result.

THEOREM 4. *Let $a, b \geq 1$ and define*

$$c(a,b) \equiv c(a,b,w) := \max_{k \geq a+b} \left( \frac{w(k-a) + w(a)}{w(k-a)} \right)$$

*and*

$$\bar{c}(a,b) \equiv c(a,b,w) := \max_{k \geq b, k \text{ even}} \left( \frac{w(k) + w(a)}{w(k)} \right) \vee \max_{k \geq a, k \text{ odd}} \left( \frac{w(k) + w(b)}{w(k)} \right).$$

*Then, for any $k \geq a + b$ and $\ell \geq a \wedge b$, we have*

$$\frac{w(a)}{c(a,b)} P^{a,b}(Z_\infty = 1) \frac{w(k-\ell) + w(\ell)}{w(k-\ell)w(\ell)} \leq P^{a,b}(Z_k = \ell) \leq w(a) \frac{w(k-\ell) + w(\ell)}{w(k-\ell)w(\ell)}$$

*and, assuming that $a$ is odd while $b$ is even,*

$$\frac{w(a) \wedge w(b)}{\bar{c}(a,b)} \prod_{j \in \{a,b\}} \bar{P}^{a,b}(Z_\infty = j)\frac{w(k-\ell) + w(\ell)}{w(k-\ell)w(\ell)}$$

$$\leq \bar{P}^{a,b}(Z_k = \ell) \leq (w(a) \vee w(b))\frac{w(k-\ell) + w(\ell)}{w(k-\ell)w(\ell)}.$$

These pre-asymptotic estimates will be useful in further analysis.

Before continuing, we note that a direct consequence is the following result, already mentioned in the Introduction.

COROLLARY 5.  *Assuming (A0), under any of the laws from Theorem 4, there exist $c, C \in (0, \infty)$, depending on the choice of the law, such that*

$$P(Z_\infty = \ell) \in \left(\frac{c}{w(\ell)}, \frac{C}{w(\ell)}\right).$$

2.2. *The "time-line" construction.*  Here, we briefly recall the construction of the edge-reinforced walk using independent families of exponentials; see [1, 11] or [5]. In the current work, we will use it mainly in the context of trees. To simplify the notation, we focus on two cases, where $\mathcal{G}$ is either a two-edge graph or a "star" with $m$ fingers. The reader can easily handle the general case.

First, assume that $\mathcal{G}$ contains two vertices, 0 and 1, and two edges, $e_G$ and $e_R$, connecting them. Fix initial weights $\ell_0^{e_G} = a$ and $\ell_0^{e_R} = b$, and let $I_{t_0} = 0$, where $t_0 = a + b$. Note that the corresponding edge-reinforced random walk has the law $P^{a,b}$.

For each $k \geq 1$, let $E_k^G$, $E_k^R$ be exponential [rate $w(k)$] random variables and let $\{E_k^G, k \geq 1\}$, $\{E_k^R, k \geq 1\}$ be two independent families of independent random variables. Let

(10)           $$T_\infty^G := \sum_{k \geq 0} E_{a+k}^G, \qquad T_\infty^R := \sum_{k \geq 0} E_{b+k}^R.$$

Note that $T_\infty^G$ and $T_\infty^R$ are independent and finite [due to (A0)] almost surely. In Figure 1, intervals between subsequent dots have length $E_{a+k}^G$ or $E_{b+k}^R$, corresponding to the edge and to the index of the chronological order $k$, and the limits $T_\infty^G$, $T_\infty^R$ are also indicated.

One can construct a realization of the edge-reinforced random walk on $\mathcal{G}$ from the above data, or (informally) from the figure, as follows.

Find the minimum of $E_a^G$ and $E_b^R$ by "simultaneously erasing at rate 1 in the chronological direction" the time-lines corresponding to both edges until the first dot is encountered. In the figure, this happens to be the first dot on the time-line corresponding to edge $e_G$, that is, $E_a^G < E_b^R$. Thus, the particle moves from 0 to 1, traversing the edge $e_G$ in the first step. Note that, due to the properties of exponentials, the probability of this move is exactly $w(a)/(w(a) + w(b))$. Continue by simultaneous erasing (the previously unerased parts of) time-lines corresponding both edges until the next dot is encountered. In the figure, it appears on the time-line of $e_R$. Hence, the particle traverses the edge $e_R$ in the second step to go back from 1 to 0. Due to the memoryless properties of exponentials, the (residual) length of the interval until the first dot on the time-line of $e_R$ is again distributed as an exponential [rate $w(b)$] random variable, independent of all other data. Therefore, the probability of this transition [namely, $w(b)/(w(b) + w(a+1))$] again matches that of the edge-reinforced random walk. Continue the above procedure of simultaneous erasure of the time-lines. In this way, the steps of the corresponding edge-reinforced random walk are generated inductively.

As its byproduct, a "continuized" version of the edge-reinforced random walk arises: here, the particle makes the jumps at exactly the times when the dots are encountered. If we denote the position of the particle (in the new process) at time $s$ by $\widetilde{I}(s)$, and if $\tau_0 = 0$ and $0 < \tau_1 < \tau_2 < \cdots$ are the subsequent jump times of the particle, then the discrete-time edge-reinforced random walk constructed above and its continuized version are coupled as follows:

$$I_k \equiv \widetilde{I}(\tau_k), \qquad k \geq 0.$$

One typically says that $I$ is the *skeleton process* of $\widetilde{I}$.

Now, let $\mathcal{G}$ be a labeled tree with the central vertex 0 which is connected via edge $e_i$ to each leaf vertex $i$, $i = 1, \ldots, m$. We call such $\mathcal{G}$ a *star with m fingers*. Fix initial weights $\ell_0^{e_i} = \ell_0^i$ and let $I_{t_0} = 0$, where $t_0 = \sum_{i=1}^m \ell_0^i$. In particular, note that if $m = 2$ and $\ell_0^1 = a$, $\ell_0^2 = b$, then the corresponding edge-reinforced random walk has the law $\bar{P}^{a,b}$.

Similarly to the previous construction, for each $i = 0, \ldots, m$ and $k \geq 1$, let $E_k^i$ be an exponential [rate $w(k)$] random variable and let $\{E_k^i, i = 0, \ldots, m, k \geq 1\}$ be a family of independent random variables. Define

$$T_\infty^i := \sum_{k \geq 0} E_{\ell_0^i + 2k}^i, \qquad i = 0, \ldots, m.$$

The multiple "2" in the subscript comes from the fact that the particle traverses each edge twice before coming back to the central vertex. Again, note that the above $m$ random variables are continuous, independent and finite almost surely. In Figure 2, intervals between subsequent dots have length $E^i_{\ell_0+2k}$ for the corresponding $i$ and $k$.

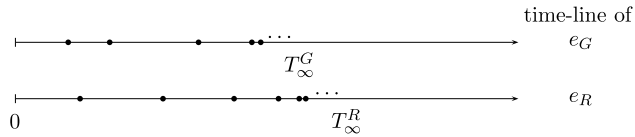

Fig. 1.  *The time-lines of a two-vertex graph.*

One constructs a realization of the corresponding edge-reinforced random walk from the above data analogously to the two-edge setting, the only difference being that now, in every second step, when not at the central vertex 0, the particle jumps almost surely back to 0. From the above figure, one can read off the first four steps of the walk as $I_{t_0+1} = 1$, $I_{t_0+1} = 0$, $I_{t_0+1} = 2$ and $I_{t_0+1} = 0$. The reader will quickly verify that, due to the properties of exponentials, the probability of transitions match those of the edge-reinforced random walk.

Again, a continuized version of the edge-reinforced random walk emerges, where there are various possibilities to account for the "singular" behavior of the walk at the leaves of $\mathcal{G}$. For example, one could use the random variables $E^i_{\ell_0+2k+1}$ (that did not play any role in the construction of the walk) as subsequent waiting times at the leaf $i$ for each $i = 1, \ldots, m$.

PROOF OF (5). We concentrate on the case $P = P^{1,1}$ and show that the limit in (5) equals $2 \int_0^\infty f(x)^2 \, dx$, where $f$ is the (continuous) density of $T^G_\infty = \sum_{k=1}^\infty E^G_k$. Let $S^R_\ell := \sum_{k=1}^{\ell-1} E^R_k$ and $f_\ell$ be the density of $S^R_\ell$. Then, since $P(Z_\infty = \ell) = P(G_\infty = \ell) + P(R_\infty = \ell) = 2P(R_\infty = \ell)$, by symmetry, we have

$$
\begin{aligned}
w(\ell)P(Z_\infty = \ell) &= 2w(\ell)P(S^R_\ell < T^G_\infty < S^R_{\ell+1}) \\
&= 2w(\ell)E(1_{\{S^R_\ell < T^G_\infty\}}P(S^R_{\ell+1} > T^G_\infty | S^R_\ell, T^G_\infty)) \\
&= 2w(\ell)E(1_{\{S^R_\ell < T^G_\infty\}}e^{-w(\ell)(T^G_\infty - S^R_\ell)}) \\
&= 2\int_0^\infty dt\, f(t) \int_0^t ds\, f_\ell(s)w(\ell)e^{-w(\ell)(t-s)},
\end{aligned}
$$
(11)

where the first identity is clear from the graphical construction above, the second is a simple conditioning relation, the third uses the fact that $S^R_{\ell+1} -$ $S_\ell^R = E_\ell^R$ is exponential [rate $w(\ell)$], independent of the $\sigma$-field generated by $S_\ell^R, T_\infty^G$, and the last is the same expression written in terms of densities of $S_\ell^R$ and $T_\infty^G$.

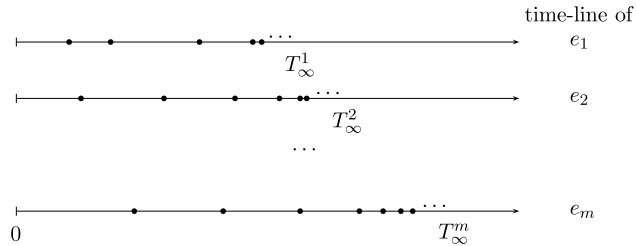

FIG. 2.    *The time-lines of a star with $m$ fingers.*

In order to prove that the integral in (11) converges to $\int_0^\infty f(x)^2\,dx$, it suffices to show that, as $\ell \to \infty$,

$$(12) \qquad \int_0^t ds\, f_\ell(s) w(\ell) e^{-w(\ell)(t-s)} \to f(t)$$

and that the left-hand side above is uniformly bounded in $\ell$ and in $t$. In fact, $f_\ell(s)$ is bounded by a fixed constant in both $\ell$ and $s$, as a convolution of an exponential and another density. Moreover,

$$|f_\ell(s_1) - f_\ell(s_2)| \leq \int |f_1(s_1 - u) - f_1(s_2 - u)| g_\ell(u)\,du,$$

where $f_1$ is the (exponential) density of $S_1^R$ and $g_\ell$ is the density of $S_\ell^R - S_1^R$. We conclude that $(f_\ell, \ell \geq 1)$ is a uniformly continuous family of functions.

The convergence in (12) is now not difficult to show, by writing $f_\ell(s) = (f_\ell(s) - f_\ell(t)) + (f_\ell(t) - f(t)) + f(t)$ and using the fact that $S_\ell^R \nearrow T_\infty^R \stackrel{d}{=} T_\infty^G$, as well as $w(\ell) \to \infty$, so that the integral concentrates around $t$ as $\ell \to \infty$.

Note that one can modify the above proof to show an analogous statement under any law $P^{a,b}$, $a, b \geq 1$, by instead using $S_\ell^G := \sum_{k=a}^{\ell-1} E_k^G$, $S_\ell^R := \sum_{k=b}^{\ell-1} E_k^R$, $T_\infty^G := \sum_{k \geq 0} E_{a+k}^G$ and $T_\infty^R := \sum_{k \geq 0} E_{b+k}^R$. Namely, one can easily verify, using the above method, that

$$\lim_{\ell \to \infty} w(\ell) P(Z_\infty = \ell) = 2 \int f_{T_\infty^G}(t) f_{T_\infty^R}(t)\,dt,$$

where $f_{T_\infty^G}$ and $f_{T_\infty^R}$ are densities of $T_\infty^G$ and $T_\infty^R$, respectively.  □

2.3. *Time of attraction.* Next, consider the time of attraction

$$T := \min\{k \geq 1 : Z_l = Z_{l+1}, \text{for all } l \geq k\}.$$

Note that $\{T = k+1\}$ is a disjoint union of $A_k^s$ and $A_k^w$, where

$A_k^s = \{$the less (or equally) traversed edge is chosen

at time $k$ and the remaining edge is chosen at all later times$\}$,

$A_k^w = \{$the more traversed edge is chosen at time $k$

and the remaining edge is chosen at all later times$\}$.

Given $\{Z_k = \ell\}$ for some $\ell < k/2$, the event $A_k^s$ happens with probability

$$(13) \qquad \frac{w(\ell)}{w(\ell) + w(k - \ell)} \prod_{j=k}^\infty \frac{w(j - \ell)}{w(j - \ell) + w(\ell + 1)},$$

while $A_k^w$ happens with probability

$$(14) \qquad \frac{w(k-\ell)}{w(k-\ell)+w(\ell)} \prod_{j=\ell}^{\infty} \frac{w(j)}{w(j)+w(k-\ell+1)} e^{-w(k)c}.$$

Finally, if $2\ell = k$, then $P(A_k^s | Z_k = \ell)$ is an expression analogous to (13),

$$\prod_{j=k}^{\infty} \frac{w(j-k/2)}{w(j-k/2)+w(k/2+1)}.$$

It will be useful to abbreviate

$$(15) \qquad \begin{aligned} W_k(\ell) &:= \prod_{j=0}^{\infty} \frac{w(k+j-\ell)}{w(k+j-\ell)+w(\ell)}, \\ \overline{W}_k(\ell) &:= \prod_{j=0}^{\infty} \frac{w(k+2j-\ell)}{w(k+2j-\ell)+w(\ell)}. \end{aligned}$$

Note that if $\ell < k/2$, then

$$(16) \qquad W_k(\ell) = P(Z_\infty = \ell | Z_k = \ell).$$

The identities (13)–(15) then yield

$$(17) \qquad P(A_k^s | Z_k = \ell) = \frac{w(\ell)}{w(\ell)+w(k-\ell)} W_{k+1}(\ell+1), \qquad \ell < k/2,$$

$$(18) \qquad P(A_k^s | Z_k = \ell) = W_{k+1}(k/2+1), \qquad k = 2\ell,$$

$$(19) \qquad P(A_k^w | Z_k = \ell) = \frac{w(k-\ell)}{w(k-\ell)+w(\ell)} W_{k+1}(k+1-\ell)$$

and similar identities hold under the laws $P^{a,b}$ and $\bar{P}^{a,b}$ (with $\overline{W}$ used in place of $W$) for $a, b \geq 1$.

We have, as discussed above,

$$\begin{aligned} P(T = k+1) = \sum_{j=1}^{\lceil k/2 \rceil - 1} P(Z_k = j)(P(A_k^s | Z_k = j) + P(A_k^w | Z_k = j)) \\ + P(Z_k = k/2) P(A_k^s | Z_k = k/2). \end{aligned}$$

Now, (17)–(19), together with Theorem 4, imply the following asymptotic formula.

LEMMA 6. *Under the law $P^{a,b}$, we have*

$$P(T = k+1) \asymp \sum_{\ell = a \wedge b}^{k/2} \left[ \frac{1}{w(k-\ell)} W_{k+1}(\ell+1) + \frac{1}{w(\ell)} W_{k+1}(k+1-\ell) \right].$$

*Similarly, under $\bar{P}^{a,b}$, where a is odd and b is even, and when k is odd (since the initial time is $a+b$ and the initial position is 0),*

$$P(T = k+2) \asymp \sum_{\ell=a,\text{odd}}^{k-b} \frac{1}{w(k-\ell)} \overline{W}_{k+2}(\ell+2) + \sum_{\ell=b,\text{even}}^{k-a} \frac{1}{w(k-\ell)} \overline{W}_{k+2}(\ell+2).$$

From now on, we will also assume that

(A1)                    $w(k)$ is nondecreasing in $k$.

This will be useful for future estimates since we then have the following.

LEMMA 7.  (a) *For each $k \geq 1$, $W_k(\cdot)$ is a nonincreasing function on the interval $[1, k/2]$.*

(b) *For each $\ell$, $W.(\ell)$ is a nondecreasing function on the interval $[2\ell, \infty)$.*

(c) *For any $k, \ell \geq 1$ such that $\ell \leq k/2$, we have $W_k(\ell) \geq W_{k+1}(\ell+1)$.*

(d) *For each fixed $\ell$, we have $\lim_{k \to \infty} W_k(\ell) \to 1$.*

PROOF.  Note that

(20)
$$W_k(\ell) = \prod_{j=0}^{\infty} \left(1 + \frac{w(\ell)}{w(k+j-\ell)}\right)^{-1},$$

$$W_{k+1}(\ell+1) = \prod_{j=0}^{\infty} \left(1 + \frac{w(\ell+1)}{w(k+j-\ell)}\right)^{-1}.$$

(a) We need to show that $W_k(\ell) \geq W_k(\ell+1)$ for all $\ell \in [1, k/2]$. By assumption (A1), we have

$$w(k+j-\ell)w(\ell+1) \geq w(k+j-\ell-1)w(\ell) \qquad \text{for all } j \geq 0,$$

from which we get that, for all $j \geq 0$,

$$\left(1 + \frac{w(\ell)}{w(k+j-\ell)}\right)^{-1} \geq \left(1 + \frac{w(\ell+1)}{w(k+j-\ell-1)}\right)^{-1}.$$

(b) Here, we need to show that $W_k(\ell) \leq W_{k+1}(\ell)$ for all $k \geq 1$. By assumption (A1), we have $w(k+j-\ell) \leq w(k+1+j-\ell)$ for all $j \geq 0$, which, together with (20), directly implies the above inequality.

(c) Again by assumption (A1), we have $w(\ell) \leq w(\ell+1)$ for all $\ell \geq 1$ so that representation (20) implies the claim.

(d) This is an easy consequence of (6) using probabilistic interpretation (16) [equivalently, one can use the algebraic definition and (A0)].  □

COROLLARY 8.   *Assume (A1). Then, for $\ell < k/2$, we have*

$$P(A_k^w | Z_k = \ell) \le 2P(A_k^s | Z_k = \ell).$$

PROOF.   Use (17) and (19). Note that

$$w(k-\ell) + w(\ell+1) \le 2(w(k+1-\ell) + w(\ell))$$

for all $k$, due to assumption (A1), and that for each $i \ge 1$, the $i$th term

$$\frac{w(\ell+i)}{w(\ell+i) + w(k-\ell+1)}$$

in the infinite product of (14) is bounded above by the $i$th term

$$\frac{w(k+i-\ell)}{w(\ell+1) + w(k+i-\ell)}$$

in the infinite product of (13), again since (A1) holds.   □

Therefore, to obtain asymptotic (in the sense of relation $\asymp$) upper and lower bounds on $P(T = k+1)$, it suffices to study only

$$P(T = k+1, A_k^s) = \sum_{\ell=1}^{k/2} P(Z_k = \ell) P(A_k^s | Z_k = \ell)$$

$$\asymp \sum_{\ell=1}^{k/2} \frac{1}{w(k-\ell)} W_{k+1}(\ell+1),$$

implying the following result.

THEOREM 9.   *If (A0) and (A1) hold, and if $P$ is $P^{a,b}$ for some fixed $a, b \ge 1$, then*

$$(21) \qquad P(T = k+1) \asymp \sum_{\ell=1}^{k/2} \frac{1}{w(k-\ell)} W_{k+1}(\ell+1)$$

*and*

$$(22) \qquad P(T > n) \asymp \sum_{k=n+1}^{\infty} \sum_{\ell=1}^{k/2} \frac{1}{w(k-\ell)} W_{k+1}(\ell+1).$$

Analogous statements are valid if $P$ is $\bar{P}^{a,b}$, where $W$ needs to be replaced by $\overline{W}$. Namely, assume that $a$ and $b$ are of opposite parity. Then, one shows in a similar fashion to Corollary 8 that for $k \ge a+b$ odd, we have $\bar{P}^{a,b}(\bar{A}_k^w | Z_k = \ell) \le 2\bar{P}^{a,b}(\bar{A}_k^s | Z_k = \ell)$, where $\bar{A}_k^s$ is the event on which the walk traverses the less (or equally) traversed edge at time $k$ and the remaining (stronger) edge from time $k+2$ onwards, and where $\bar{A}_k^w = \{T = k+2\} \setminus \bar{A}_k^s$. Combining this with Lemma 6, we obtain

$$\bar{P}(T = k+2) \asymp \sum_{\ell=1}^{k/2} \frac{1}{w(k-\ell)} \overline{W}_{k+2}(\ell+2)$$

and

$$\bar{P}(T > n) \asymp \sum_{\substack{k=n+1, \\ k \text{ odd}}}^{\infty} \sum_{\ell=1}^{k/2} \frac{1}{w(k-\ell)} \overline{W}_{k+2}(\ell+2).$$

2.4. *Examples.* Let $a, b \geq 1$ be fixed integers and denote by $P$ either of the laws $P^{a,b}$ or $\bar{P}^{a,b}$.

THEOREM 10. *Suppose $\rho > 1$, let $\rho' := (\rho-1)/\rho$, $\alpha \in (0, \infty)$ and let $\epsilon > 0$ be arbitrarily fixed.*

(a) *If $w(k) = k^\rho$, $k \geq 1$, then there exist finite positive $c_1(\rho), c_2(\rho)$ such that for all $k \geq 1$,*

$$\frac{c_1(\rho)}{k^{\rho-\rho'}} \leq P(T = k+1) \leq \frac{c_2(\rho)}{k^{\rho-\rho'}}.$$

(b) *If $w(k) = k^\rho \log^\alpha k$, $k \geq 1$, $\alpha > 1$, then there exist finite positive $c_1(\alpha, \rho)$, $c_2(\alpha, \rho)$ such that for all $k \geq 1$,*

$$\frac{c_1(\alpha, \rho)}{k^{\rho-\rho'} \log^\alpha k} \leq P(T = k+1) \leq \frac{c_2(\alpha, \rho)}{k^{\rho-\rho'} \log^\alpha k}.$$

(c) *If $w(k) = ke^{\log^\alpha k}$, $k \geq 1$, $0 < \alpha < 1$, then there exist finite positive $c_1(\alpha), c_2(\alpha)$ such that for all $k \geq 1$ and all $\epsilon > 0$,*

$$\frac{c_1(\alpha)}{k(\log k)^{1-\alpha} e^{\log^{\alpha^2} k - \beta \log^{2\alpha^2 - \alpha} k}} \leq P(T = k+1)$$

$$\leq \frac{c_2(\alpha)}{k(\log k)^{1-\alpha} e^{\log^{\alpha^2} k - (\alpha+\epsilon) \log^{2\alpha^2 - \alpha} k}},$$

*where $\beta = \alpha$ if $0 < \alpha \leq 2/3$ and $\beta = \alpha - \epsilon$ otherwise.*

*In particular, if $0 < \alpha \leq 1/2$, then $\exp\{-\log^{2\alpha^2 - \alpha} k\} \asymp 1$, so*

$$P(T = k+1) \asymp \frac{c_2(\alpha)}{k(\log k)^{1-\alpha} e^{\log^{\alpha^2} k}}.$$

(d) *If $w(k) = e^{k^\epsilon}$ for some $0 < \epsilon \le 1$, then there exist finite positive $c_1(\epsilon), c_2(\epsilon)$ such that for all $k \ge 1$,*

$$\frac{c_1(\epsilon)k^{1-\epsilon}}{e^{(k-k^*)^\epsilon}} \le P(T = k+1) \le \frac{c_2(\epsilon)k^{1-\epsilon}}{e^{(k-k^*)^\epsilon}},$$

*where $k^* = \frac{k}{2} - \frac{(1-\epsilon)}{\epsilon 2^{2-\epsilon}} k^{1-\epsilon} \log k$.*

(e) *If $w(k) = e^{\alpha k}$, then there exist finite positive $c_1(\alpha), c_2(\alpha)$ such that*

$$c_1(\alpha)e^{-\alpha k/2} \le P(T = k+1) \le c_2(\alpha)e^{-\alpha k/2}.$$

REMARK. All of the constants of the form $c_1(\cdot), c_2(\cdot)$ featuring in the statements above additionally depend on the initial weights $a$ and $b$, due to Theorems 4 and 9; for an example, see (26) below.

PROOF OF THEOREM 10. We concentrate on the case where $P = P^{a,b}$; the other case, $P = \bar{P}^{a,b}$, can be treated similarly. Without loss of generality, we assume that $\ell \le k/2$. We are going to use the inequality

$$(23) \quad e^{-w(\ell+1)\sum_{j=k}^{\infty} 1/w(j-\ell)} \le W_{k+1}(\ell+1) \le e^{-w(\ell+1)\sum_{j=k}^{\infty} 1/2w(j-\ell)},$$

which is a direct consequence of the fact that

$$e^{-x} \le (1+x)^{-1} \le e^{-x/2}, \qquad 0 \le x \le 1.$$

(a) For $w(k) = k^\rho$, (23) becomes

$$e^{-(\ell+1)^\rho \sum_{j=k}^{\infty} 1/(j-\ell)^\rho} \le W_{k+1}(\ell+1) \le e^{-(\ell+1)^\rho \sum_{j=k}^{\infty} 1/(2(j-\ell)^\rho)}.$$

Using

$$\int_{k+1}^{\infty} \frac{dx}{(x-\ell)^\rho} \le \sum_{j=k}^{\infty} \frac{1}{(j-\ell)^\rho} \le \int_{k}^{\infty} \frac{dx}{(x-\ell)^\rho},$$

we get a lower and an upper bound

$$(24) \quad \frac{1}{\rho-1}\frac{1}{(k+1-\ell)^{\rho-1}} \le \sum_{j=k}^{\infty} \frac{1}{(j-\ell)^\rho} \le \frac{1}{\rho-1}\frac{1}{(k-\ell)^{\rho-1}}.$$

Therefore, we have, for $1 \le \ell \le k/2$,

$$e^{-(2^{\rho-1}/(\rho-1))((\ell+1)^\rho/k^{\rho-1})} \le e^{-(1/(\rho-1))((\ell+1)^\rho/(k-\ell)^{\rho-1})} \le W_{k+1}(\ell+1)$$

$$(25) \qquad\qquad \le e^{-(1/(2\rho-2))((\ell+1)^\rho/(k+1-\ell)^{\rho-1})}$$

$$\qquad\qquad \le e^{-(1/(2\rho-2))(\ell^\rho/k^{\rho-1})}.$$

Now, Theorem 9 implies that

$$(26) \qquad P(T = k+1) \geq c(a,b) \sum_{\ell=1}^{k^{\rho'}} \frac{1}{(k-\ell)^\rho} e^{-(2^{\rho-1}/(\rho-1))((\ell+1)^\rho/k^{\rho-1})}$$

$$\geq c(a,b) e^{-(2^{2\rho-1}/(\rho-1))} \sum_{\ell=1}^{k^{\rho'}} \frac{1}{(k-\ell)^\rho}$$

$$(27) \qquad \geq c(a,b) e^{-(2^{2\rho-1}/(\rho-1))} k^{\rho'-\rho},$$

where $c(a,b)$ is a finite positive constant. For the upper bound, again use Theorem 9 and (25). Now, note that, since $k \geq 2\ell$,

$$(28) \qquad \sum_{\ell=k^{\rho'}}^{k/2} \frac{1}{(k-\ell)^\rho} e^{-(1/(2\rho-2))(\ell^\rho/k^{\rho-1})} \leq \left(\frac{2}{k}\right)^\rho \sum_{\ell=k^{\rho'}}^{k/2} e^{-(1/(2\rho-2))(\ell^\rho/k^{\rho-1})}.$$

To bound the last term above, we split the interval $[k^{\rho'}, k/2]$ into subintervals of equal width $k^{\rho'}$, with the last subinterval possibly having smaller width $k/2 - \lfloor k^{1-\rho'}/2 \rfloor k^{\rho'}$. To abbreviate, we define $a_k := \lfloor k^{1-\rho'}/2 \rfloor$. We then have

$$\sum_{\ell=k^{\rho'}}^{k/2} e^{-(1/(2\rho-2))(\ell^\rho/k^{\rho-1})}$$

$$\leq k^{\rho'} e^{-1/(2\rho-2)} + k^{\rho'} e^{-2^\rho/(2\rho-2)} + \cdots + k^{\rho'} e^{-a_k{}^\rho/(2\rho-2)}$$

$$\leq k^{\rho'} [e^{-1/(2\rho-2)} + e^{-2/(2\rho-2)} + \cdots + e^{-a_k/(2\rho-2)}]$$

$$\leq \frac{k^{\rho'}}{1 - e^{-1/(2\rho-2)}},$$

which, together with (28), completes the proof of part (a).

(b) Using (23) and the fact that

$$\sum_{j=k}^{\infty} \frac{1}{(j-\ell)^\rho \log^\alpha(j-\ell)}$$

is up to a constant multiple of order $1/((k-\ell)^{\rho-1} \log^\alpha(k-\ell))$, one obtains

$$e^{-c_1(\alpha,\rho)((\ell+1)^\rho/(k-\ell)^{\rho-1})(\log^\alpha(\ell+1)/\log^\alpha(k-\ell))}$$

$$(29) \qquad \leq W_{k+1}(\ell+1)$$

$$\leq e^{-c_2(\alpha,\rho)((\ell+1)^\rho/k^{\rho-1})(\log^\alpha(\ell+1)/\log^\alpha k)},$$

where $c_1(\alpha,\rho), c_2(\alpha,\rho)$ are finite positive constants. As in the case (a), one gets the lower bound by evaluating the order of

$$(30) \qquad \sum_{\ell=\xi}^{\zeta} \frac{1}{w(k-\ell)} W_{k+1}(\ell+1)$$

for $\xi = 1$, $\zeta = k^{\rho'}$ and the upper bound by evaluating the order of the sum for $\xi = 1$, $\zeta = k^{\rho'}$ and for $\xi = k^{\rho'}$, $\zeta = k/2$, separately.

(c) We have

$$\int_{k+1}^{\infty} \frac{dx}{(x-\ell)e^{\log^{\alpha}(x-\ell)}} \leq \sum_{j=k}^{\infty} \frac{1}{(j-\ell)e^{\log^{\alpha}(j-\ell)}} \leq \int_{k}^{\infty} \frac{dx}{(x-\ell)e^{\log^{\alpha}(x-\ell)}}.$$

Since

$$\int_{k}^{\infty} \frac{dx}{(x-\ell)e^{\log^{\alpha}(x-\ell)}}$$
$$= \frac{\log^{1-\alpha}(k-\ell)}{\alpha e^{\log^{\alpha}(k-\ell)}} + \frac{1-\alpha}{\alpha} \int_{k}^{\infty} \frac{dx}{(x-\ell)\log^{\alpha}(x-\ell)e^{\log^{\alpha}(x-\ell)}}$$
$$= \frac{\log^{1-\alpha}(k-\ell)}{\alpha e^{\log^{\alpha}(k-\ell)}}(1 + o_{k-\ell}(1)),$$

where $o_k(1) \to 0$ as $k \to \infty$, using the fact that $1 \leq \ell \leq \frac{k}{2}$ above, we get the inequality

$$c_1(\alpha)\frac{(\log k)^{1-\alpha}}{e^{\log^{\alpha} k}} \leq \sum_{j=k}^{\infty} \frac{1}{(j-\ell)e^{\log^{\alpha}(j-\ell)}} \leq c_2(\alpha)\frac{(\log k)^{1-\alpha}}{e^{\log^{\alpha} k}}$$

for some $c_1(\alpha), c_2(\alpha) \in (0,\infty)$. By now applying (23), we obtain

$$e^{-c_1(\alpha)(\ell+1)e^{\log^{\alpha}(\ell+1)}\log^{1-\alpha} k/(e^{\log^{\alpha} k})}$$
$$\leq W_{k+1}(\ell+1) \leq e^{-c_2(\alpha)(\ell+1)e^{\log^{\alpha}(\ell+1)}\log^{1-\alpha} k/(e^{\log^{\alpha} k})}.$$

As in parts (a) and (b), we find a convenient breaking point and approximate the sums (30) separately. We take

$$\zeta = e^{\log^{\alpha} k - \log^{\alpha^2} k + \beta \log^{2\alpha^2 - \alpha} k}/\log^{1-\alpha} k,$$

where $\beta = \alpha$ for $0 < \alpha \leq 2/3$ and $\beta = \alpha - \epsilon$ otherwise.

To verify the lower bound, we need to show that we can bound

$$\exp\left\{-c_1(\alpha)(\zeta+1)e^{\log^{\alpha}(\zeta+1)}\frac{\log^{1-\alpha} k}{e^{\log^{\alpha} k}}\right\}$$

from below by a positive constant. Hence, we estimate (the constant $c$ below is finite and positive, and possibly changes from line to line)

$$\frac{(\zeta+1)\log^{1-\alpha} k}{\exp(\log^{\alpha} k)}\exp(\log^{\alpha}(\zeta+1))$$

(31)
$$\leq c\exp[\beta \log^{2\alpha^2-\alpha} k - \log^{\alpha^2} k]$$
$$\times \exp[(\log^{\alpha} k - \log^{\alpha^2} k + \beta \log^{2\alpha^2-\alpha} k - (1-\alpha)\log\log k)^{\alpha}]$$

$$\leq c \exp[\beta \log^{2\alpha^2-\alpha} k - \log^{\alpha^2} k$$
$$+ \log^{\alpha^2} k (1 - \log^{\alpha^2-\alpha} k + \beta \log^{2\alpha^2-2\alpha} k)^\alpha]$$

$$\leq c \exp[\beta \log^{2\alpha^2-\alpha} k - \log^{\alpha^2} k$$
$$+ \log^{\alpha^2} k (1 - \alpha \log^{\alpha^2-\alpha} k + \beta \alpha \log^{2\alpha^2-2\alpha} k)]$$

$$\tag{32} = c \exp(\beta \log^{2\alpha^2-\alpha} k - \alpha \log^{2\alpha^2-\alpha} k + \beta \alpha \log^{3\alpha^2-2\alpha} k),$$

where $\beta$ is chosen as above and where, for the third inequality, we use the fact that $(1-x)^\epsilon \leq 1 - \epsilon x$ for $0 \leq x < 1$ and $0 \leq \epsilon \leq 1$. Note, now, that if $\alpha \leq 2/3$, then $\log^{3\alpha^2-2\alpha} k$ is bounded from above by a constant, so $\beta = \alpha$ is sufficient to bound the expression (32) from above by a constant. If $\alpha > 2/3$, then $\log^{3\alpha^2-2\alpha} k \to \infty$ and $3\alpha^2 - 2\alpha < 2\alpha^2 - \alpha$, so taking $\beta = \alpha - \epsilon$ for any $\epsilon > 0$ will suffice to bound (32) by a constant. By Lemma 7(a), $W_{k+1}(\ell+1) \geq W_{k+1}(\zeta+1) \geq e^{-c'}$, for $\ell \leq \zeta$, so

$$\tag{33} \sum_{\ell=1}^{\zeta} \frac{1}{w(k-\ell)} W_{k+1}(\ell+1) \geq e^{-c'} \frac{\zeta}{w(k)} \geq \frac{c_1(\alpha)}{k \log^{1-\alpha} k e^{\log^{\alpha^2} k - \beta \log^{2\alpha^2-\alpha}}},$$

which proves the lower bound.

For the upper bound, break the summation at the point

$$\zeta := e^{\log^\alpha k - \log^{\alpha^2} k + (\alpha+\epsilon) \log^{2\alpha^2-\alpha} k} / \log^{1-\alpha} k.$$

We have, since $W_{k+1}(\ell+1) \leq 1$,

$$\sum_{\ell=1}^{\zeta} \frac{1}{w(k-\ell)} W_{k+1}(\ell+1) \leq \sum_{\ell=1}^{\zeta} \frac{1}{w(k-\ell)} \leq \frac{\zeta}{w(k-\zeta)}$$

and since $w(k-\zeta) \asymp w(k)$, we can bound this term from above by a term of the order stated in the formulation of part (c). Hence, it suffices to bound $\sum_{\ell=\zeta}^{k/2} W_{k+1}(\ell+1)/w(k-\ell)$ as follows:

$$\sum_{\ell=\zeta}^{k/2} \frac{1}{w(k-\ell)} W_{k+1}(\ell+1) \leq \frac{1}{w(k/2)} \sum_{\ell=\zeta}^{k/2} W_{k+1}(\ell+1)$$

$$\leq \frac{c_2}{k e^{\log^\alpha k}} \int_\zeta^\infty \frac{dx}{e^{c_2(\alpha)x e^{\log^\alpha x} \log^{1-\alpha} k / e^{\log^\alpha k}}}$$

$$\tag{34} = \frac{c_2}{k e^{\log^\alpha k}} \cdot \frac{e^{\log^\alpha k}}{\log^{1-\alpha} k (1 + \alpha \log^{\alpha-1}(\zeta)) e^{\log^\alpha(\zeta)}}$$

$$\cdot \frac{1}{e^{c_2(\alpha) \cdot \zeta \cdot e^{\log^\alpha(\zeta)} \log^{1-\alpha} k / e^{\log^\alpha k}}} (1 + o_\zeta(1)),$$

where the last term is obtained by integration by parts.

First, estimate

$$\frac{\zeta \cdot \log^{1-\alpha} k}{\exp(\log^{\alpha} k)} \exp(\log^{\alpha}(\zeta))$$

$$= \exp[(\alpha + \epsilon) \log^{2\alpha^2 - \alpha} k - \log^{\alpha^2} k]$$

$$\times \exp[(\log^{\alpha} k - \log^{\alpha^2} k + (\alpha + \epsilon) \log^{2\alpha^2 - \alpha} k - (1 - \alpha) \log \log k)^{\alpha}]$$

$$\geq c \exp[(\alpha + \epsilon) \log^{2\alpha^2 - \alpha} k - \log^{\alpha^2} k + \log^{\alpha^2} k (1 - \log^{2\alpha^2 - \alpha} k)^{\alpha}]$$

$$\geq c \exp[(\alpha + \epsilon) \log^{2\alpha^2 - \alpha} k - \log^{\alpha^2} k$$

$$+ \log^{\alpha^2} k (1 - (\alpha + \epsilon) \log^{\alpha^2 - \alpha} k)] \geq c',$$

where, in the last inequality above, we used

$$(35) \qquad\qquad (1 - x)^{\alpha} \geq 1 - (\alpha + \epsilon)x, \qquad \epsilon > 0, x \to 0.$$

Therefore, we can bound the multiple (34) from above by a constant $c$. Furthermore,

$$\frac{c_2}{k e^{\log^{\alpha} k}} \cdot \frac{e^{\log^{\alpha} k}}{\log^{1-\alpha} k (1 + \alpha \log^{\alpha-1}(\zeta)) e^{\log^{\alpha}(\zeta)}}$$

$$\leq \frac{c_2 \exp\{-(\log^{\alpha} k - \log^{\alpha^2} k + (\alpha + \epsilon) \log^{2\alpha^2 - \alpha} k - (1 - \alpha) \log \log k)^{\alpha}\}}{k (\log k)^{1-\alpha}}$$

$$\leq \frac{c_2 \exp\{-(\log^{\alpha} k - \log^{\alpha^2} k)^{\alpha}\}}{k (\log k)^{1-\alpha}}$$

$$\leq \frac{c_2 \exp\{-\log^{\alpha^2} k (1 - (\alpha + \epsilon) \log^{\alpha^2 - \alpha} k)\}}{k (\log k)^{1-\alpha}},$$

where, in the last inequality above, we again used (35).

(d) We have

$$\int_{k-\ell+1}^{\infty} \frac{dx}{e^{x^{\epsilon}}} \leq \sum_{j=k-\ell}^{\infty} \frac{1}{e^{j^{\epsilon}}} \leq \int_{k-\ell}^{\infty} \frac{dx}{e^{x^{\epsilon}}}.$$

Note that

$$\int_{k-\ell}^{\infty} \frac{dx}{e^{x^{\epsilon}}} = \frac{(k-\ell)^{1-\epsilon}}{\epsilon e^{(k-\ell)^{\epsilon}}} + \frac{1-\epsilon}{\epsilon} \int_{k-\ell}^{\infty} \frac{dx}{x^{\epsilon} e^{x^{\epsilon}}} = \frac{(k-\ell)^{1-\epsilon}}{\epsilon e^{(k-\ell)^{\epsilon}}} (1 + o_{k-\ell}(1)),$$

where, again, $o_k(1) \to 0$ as $k \to \infty$. From the above formulae, we obtain the inequality

$$(36) \qquad c_2(\epsilon) \frac{(k-\ell)^{1-\epsilon}}{e^{(k-\ell)^{\epsilon}}} \leq \sum_{j=k-\ell}^{\infty} \frac{1}{e^{j^{\epsilon}}} \leq c_1(\epsilon) \frac{(k-\ell)^{1-\epsilon}}{e^{(k-\ell)^{\epsilon}}}$$

for some $c_1(\epsilon), c_2(\epsilon) \in (0, \infty)$. By now applying (23), we obtain

$$(37) \quad e^{-c_1(\epsilon)e^{\ell^\epsilon}(k-\ell)^{1-\epsilon}/e^{(k-\ell)^\epsilon}} \leq W_{k+1}(\ell+1) \leq e^{-c_2(\epsilon)e^{\ell^\epsilon}(k-\ell)^{1-\epsilon}/e^{(k-\ell)^\epsilon}}.$$

If we now take $k^* = \frac{k}{2} - \frac{(1-\epsilon)}{2\epsilon 2^{1-\epsilon}} k^{1-\epsilon} \log k$, then

$$W_{k+1}(k^*+1) \geq e^{-c_1(\epsilon)e^{k^{*\epsilon}}(k-k^*)^{1-\epsilon}/e^{(k-k^*)^\epsilon}}.$$

Note that

$$k^{*\epsilon} - (k-k^*)^\epsilon = (k-k^*)^\epsilon \left[ \left(1 - \frac{(1-\epsilon)k^{1-\epsilon}\log k}{\epsilon 2^{1-\epsilon}(k-k^*)}\right)^\epsilon - 1 \right]$$

$$\leq -(1-\epsilon)\log k + O(k^{-\epsilon/2}),$$

where the last inequality is obtained using $(1-x)^\epsilon \leq 1 - \epsilon x$ for $0 \leq x < 1$ and $0 \leq \epsilon \leq 1$, together with $k - k^* \geq k/2(1 + O(\log k/k^\epsilon))$. Therefore, $\liminf_k W_{k+1}(k^*+1) > 0$ and, by Lemma 7(a),

$$W_{k+1}(\ell+1) \geq W_{k+1}(k^*+1) \geq c(\epsilon), \qquad \ell \leq k^* \text{ for some } c(\varepsilon) > 0.$$

Therefore, recalling $c_1(\epsilon), c_2(\epsilon)$ from (36),

$$\sum_{\ell=1}^{k^*} \frac{1}{w(k-\ell)} W_{k+1}(\ell+1) \geq c(\epsilon) \sum_{\ell=1}^{k^*} \frac{1}{e^{(k-\ell)^\epsilon}} \geq c(\epsilon) \int_{k-k^*+1}^{k+1} \frac{dx}{e^{x^\epsilon}}$$

$$\geq c(\epsilon) \left( \frac{c_1(\epsilon)(k-k^*)^{1-\epsilon}}{e^{(k-k^*)^\epsilon}} - \frac{c_2(\epsilon)k^{1-\epsilon}}{e^{k^\epsilon}} \right)$$

$$\geq \frac{c(\epsilon)c_1(\epsilon)(k-k^*)^{1-\epsilon}}{e^{(k-k^*)^\epsilon}},$$

where the constant $c(\epsilon)$ may change from line to line by a positive finite multiple. This proves the lower bound.

To get the corresponding upper bound, first observe that $W_{k+1}(\ell+1) \leq 1$, so, using (36), we can simply bound

$$\sum_{\ell=1}^{k^*} \frac{1}{w(k-\ell)} W_{k+1}(\ell+1) \leq \sum_{j=k-k^*}^{\infty} \frac{1}{w(j)} \leq c_2(\epsilon) \frac{(k-k^*)^{1-\varepsilon}}{e^{(k-k^*)^{1-\varepsilon}}}$$

and we proceed to bound $\sum_{\ell=k^*}^{k/2} \frac{1}{w(k-\ell)} W_{k+1}(\ell+1)$.

Due to (36) and (37), we can write, for $\ell \in [k^*, k/2]$,

$$\exp\left\{ -c_1'(\epsilon)w(\ell) \int_{k-\ell}^{\infty} \frac{dx}{e^{x^\epsilon}} \right\} \leq W_{k+1}(\ell+1) \leq \exp\left\{ -c_2'(\epsilon)w(\ell) \int_{k-\ell}^{\infty} \frac{dx}{e^{x^\epsilon}} \right\}$$

(38)

for some $c_1'(\epsilon), c_2'(\epsilon) \in (0, \infty)$. Since $W_{k+1}(k^* + 1) \asymp 1$, as we showed in the proof of the lower bound, we have

$$\frac{1}{w(k^*)} \asymp \int_{k-k^*}^{\infty} \frac{dx}{e^{x^\epsilon}}$$

and since

$$\int_{k-\ell}^{\infty} \frac{dx}{e^{x^\epsilon}} \geq \int_{k-k^*}^{\infty} \frac{dx}{e^{x^\epsilon}}, \qquad \ell \geq k^*,$$

we conclude from (38) that

$$W_{k+1}(\ell + 1) \leq \exp\{-cw(\ell)/w(k^*)\}, \qquad \ell \in [k^*, k/2].$$

Therefore,

$$\sum_{\ell=k^*}^{k/2} \frac{1}{w(k-\ell)} W_{k+1}(\ell+1) \leq \sum_{\ell=k^*}^{k/2} \frac{e^{-cw(\ell)/w(k^*)}}{w(k-\ell)} \leq \sum_{\ell=k^*}^{k/2} \frac{\bar{c} w^2(k^*)}{w^2(\ell) w(k-\ell)},$$

where $\bar{c} \in (0, \infty)$ is such that $e^{-cx} \leq \bar{c} x^{-2}$ for all $x \geq 1$ and where we use the fact that $w(\ell) \geq w(k^*)$ for $\ell \geq k^*$. Finally, it is easy to check that $\ell \mapsto w(k-\ell)w(\ell)$, $\ell \geq k^*$, is a nondecreasing function, so

$$\frac{w(k^*)}{w(\ell)w(k-\ell)} \leq \frac{1}{w(k-k^*)}$$

and therefore

$$\sum_{\ell=k^*}^{k/2} \frac{w^2(k^*)}{w^2(\ell)w(k-\ell)} \leq \frac{w(k^*)}{w(k-k^*)} \sum_{\ell=k^*}^{k/2} \frac{1}{w(\ell)} \leq \frac{w(k^*)}{w(k-k^*)} \int_{k^*-1}^{k/2} \frac{dx}{w(x)}$$

$$\leq c_2(\epsilon) \frac{(k^*)^{1-\epsilon}}{w(k-k^*)},$$

which gives the upper bound, due to the fact that $k^* < k/2 < k - k^*$.

(e) This is a direct consequence of part (d), but its direct proof (left to an interested reader) is much easier. This fact is related to the following property: among all of the weights in (d), it is only the case of $w(k) = e^k$ where the edge-reinforced random walk gets attracted at any particular time with probability uniformly bounded away from zero. $\quad\square$

As a consequence of Theorem 10, we now have the following.

COROLLARY 11. *Suppose that $w(k)$ is as in* (a) *or* (b) *in Theorem* 10. *Then, $E(T)$ is infinite if $\rho \leq 1 + \frac{1+\sqrt{5}}{2}$ and finite if $\rho > 1 + \frac{1+\sqrt{5}}{2}$.*

**3. Analysis on general graphs.** Assume that $\mathcal{G}$ is a connected graph with $D(\mathcal{G}) < \infty$. Recall that $P^{\mathcal{G}}$ is the law of the reinforced random walk on $\mathcal{G}$.

We start with an easy lower bound in terms of the tail distribution of $T$ under the two-edge law $\bar{P}$. In fact, in the following comparison arguments, it will be convenient to instead consider the law of

$$(T)^+ \equiv T^+ := T - t_0.$$

Unless otherwise stated, in this section, we will assume that $\ell_0^e$, $e \in E(\mathcal{G})$, forms the (general) initial configuration of weights on edges such that $\ell_0^e < \infty$, $e \in E(\mathcal{G})$.

LEMMA 12. *There exist $c = c(w, \mathcal{D}(\mathcal{G})) \in (0, \infty)$ and $a, b \in \mathbb{N}$ such that*

$$P^{\mathcal{G}}(T^+ > k) \geq c\bar{P}^{a,b}(T^+ > k).$$

PROOF. Let $I_{t_0} = v \in \mathcal{G}$ be the initial position. Without loss of generality, assume that at least two edges $e$ and $f$ meet at $v$. Otherwise, at least two edges must meet at the unique neighbor of $v$ and the argument is similar. Recall that $\mathcal{G}_1$ denotes the range of the walk. Define the event

$$A_{e,f} := \{\mathcal{G}_1 \subset \text{ graph spanned by } e, f\}$$

and note that, due to (A0), event $A_{e,f}$ has positive probability for any given bounded degree graph and any fixed configuration $\ell_0^e$, $e \in E(\mathcal{G})$. At the same time,

$$\{T^+ > k\} \cap A_{e,f} \subset \{T^+ > k\}.$$

Denote by $v_e$ and $v_f$ the two vertices such that $e = \{v, v_e\}$ and $f = \{v, v_f\}$. We will verify below the existence of a positive constant $\beta$ that depends on $\mathcal{G}, w$ and the initial weights, such that for each (possibly infinite) path $v = i_0 \sim i_1 \sim \cdots$ of vertices where $i_n \in \{v, v_e, v_f\}$, $n \geq 0$, we have

$$(39) \quad P^{\mathcal{G}}(I_{t_0} = i_0, I_{t_0+1} = i_1, \ldots) \geq \beta\bar{P}^{\ell_0^e, \ell_0^f}(I_{t_0} = i_0, I_{t_0+1} = i_1, \ldots).$$

Note that $t_0$ equals $\ell_0^e, \ell_0^f$ under the law $\bar{P}^{\ell_0^e, \ell_0^f}$, but as mentioned earlier, the edge-reinforced random walk can be redefined by a time-shift to start from any fixed initial time and this does not change the probability of it taking any particular path. Clearly, (39) implies that $P^{\mathcal{G}}(B \cap A_{e,f}) \geq \beta\bar{P}^{\ell_0^e, \ell_0^f}(B)$ for any event $B$ in the $\sigma$-field generated by the walk. In particular,

$$P^{\mathcal{G}}(T^+ > k) \geq P^{\mathcal{G}}(\{T^+ > k\} \cap A_{e,f}) \geq \beta\bar{P}^{a,b}(T^+ > k),$$

as claimed.

It suffices to verify (39) for each infinite path $i_0 \sim i_1 \sim \cdots$ specified above. For $n \geq t_0$, define $x_n^e := \ell_0^e + \#\{j \leq n : \{i_{j-1}, i_j\} = e\}$ and $x_n^f := \ell_0^f + \#\{j \leq n : \{i_{j-1}, i_j\} = f\}$. The probability on the right-hand side of (39) equals

$$(40) \qquad \prod_{n=0}^{\infty} \frac{w(x_{t_0+2n+1}^{\{i_{2n}, i_{2n+1}\}})}{w(x_{t_0+2n}^e) + w(x_{t_0+2n}^f)}.$$

Define

$$c(v) := \sum_{u : u \sim v, u \neq v_e, v_f} w(\ell_0^{\{u,v\}}),$$

$$c(v_e) := \sum_{u : u \sim v_e, u \neq v} w(\ell_0^{\{u,v_e\}}), \qquad c(v_f) := \sum_{u : u \sim v_f, u \neq v} w(\ell_0^{\{u,v_f\}}).$$

The probability on the left-hand side of (39) equals

$$(41) \qquad \prod_{n=0}^{\infty} \frac{w(x_{t_0+2n}^{\{i_{2n}, i_{2n+1}\}})}{w(x_{t_0+2n}^e) + w(x_{t_0+2n}^f) + c(v)} \prod_{n=0}^{\infty} \frac{w(x_{t_0+2n+1}^{\{i_{2n+1}, i_{2n+2}\}})}{w(x_{t_0+2n+1}^{\{i_{2n+1}, i_{2n+2}\}}) + c(v_{\{i_{2n+1}, i_{2n+2}\}})},$$

where the first infinite product accounts for all the steps originating from the middle vertex $v$, while the second infinite product accounts for all the steps originating from the "boundary vertices" $v_e$ and $v_f$. Since

$$\sum_{n=0}^{\infty} \frac{c(v_{\{i_{2n+1}, i_{2n+2}\}})}{w(x_{t_0+2n+1}^{\{i_{2n+1}, i_{2n+2}\}}) + c(v_{\{i_{2n+1}, i_{2n+2}\}})} \leq 2 \sum_{n=0}^{\infty} \frac{c(v_e) \vee c(v_f)}{w(n) + (c(v_e) \wedge c(v_f))} < \infty$$

by (A0), a well-known calculus fact implies that the second product is uniformly (over infinite paths) bounded away from 0. The ratio of the first product in (40) and the probability in (42) is again uniformly bounded away from 0 since

$$(42) \qquad \sum_{n=0}^{\infty} \frac{c(v)}{w(x_{t_0+2n}^e) + w(x_{t_0+2n}^f) + c(v)} \leq 2 \sum_{n=0}^{\infty} \frac{c(v)}{w(n) + c(v)} < \infty. \qquad \square$$

Getting a corresponding upper bound on the tails of the distribution of $T$ seems more difficult. As a warm-up, we study the tree setting next, and the general finite graph and infinite graph settings, respectively, in the following subsections.

The following fact, complementary in spirit to conditioning on event $A_{e,f}$ in the proof of Lemma 12, will soon prove useful.

LEMMA 13. *Suppose that $\mathcal{G}^*$ is a finite connected subgraph of $\mathcal{G}$. Then, for each (possibly infinite) path $i_0 \sim i_1 \sim \cdots$ of vertices all contained in $\mathcal{G}^*$, we have, assuming $P^{\mathcal{G}}(I_{t_0} = i_0) = P^{\mathcal{G}^*}(I_{t_0} = i_0) = 1$,*

$$P^{\mathcal{G}}(I_{t_0} = i_0, I_{t_0+1} = i_1, \ldots) \leq P^{\mathcal{G}^*}(I_{t_0} = i_0, I_{t_0+1} = i_1, \ldots).$$

PROOF.   At each step $k$ where all the neighbors of the current position $i_k$ are contained in $\mathcal{G}^*$, the probability of the transition from $i_k$ to $i_{k+1}$ is the same under both laws $P^{\mathcal{G}}$ and $P^{\mathcal{G}^*}$.

At each step $k$ where at least one neighbor of the current position $i_k$ is an element of $V(\mathcal{G}) \setminus V(\mathcal{G}^*)$, note that the probability of the transition from $i_k$ to $i_{k+1}$ under $P^{\mathcal{G}}$ is strictly smaller than that under $P^{\mathcal{G}^*}$.   □

3.1. *Analysis on trees.*  In this subsection, we assume that $\mathcal{G}$ is a tree such that $D(\mathcal{G}) < \infty$ and we derive some upper bound estimates on the tail distribution of $T$ under $P^{\mathcal{G}}$.

First, let $\mathcal{G}$ be the star with $m$ fingers, as defined in Section 2.2. For the sake of concreteness, we assume that all of the initial weights $\ell_0^{e_i}$ are equal to 1 and that $I_{t_0} = I_m = 0$. A similar statement applies for more general initial configurations.

LEMMA 14.   $P^{\mathcal{G}}(T^+ > k) \leq \binom{m}{2} \bar{P}^{1,1}(T^+ > k/\binom{m}{2})$.

PROOF.   We will show that

$$(43) \qquad T^+ \leq \sum_{1 \leq i < j \leq m} T^{e_i, e_j, restr, +}, \qquad \text{almost surely,}$$

where $T^{e_i, e_j, restr} = T^{e_j, e_i, restr}$ is a random variable to be defined, corresponding to the pair of edges $e_i, e_j$ such that its law under $P^{\mathcal{G}}$ is the law of $T$ under $\bar{P}^{1,1}$ and where $T^{e_i, e_j, restr, +} = (T^{e_i, e_j, restr})^+$. The reason for (43) is as follows. Suppose that $f \in \{e_1, \ldots, e_m\}$ is the attracting edge for the walk. The steps away from the central vertex 0 up to time $T$ are naturally split into $K_e$ steps traversing edge $e \neq f$ (so, in total, there are $2K_e$ steps along any $e \neq f$). Up to time $T$, there are therefore $T - 2 \sum_{e \neq f} K_e$ steps across $f$. For $e \neq f$, define

$$Y^{e,f} := \text{time of the last traversal of } e$$

and

$$T_n^{e,f} := \ell_0^e + \ell_0^f + \#\text{traversals of } e \text{ or } f \text{ up to time } n,$$

$$T_n^{e,f,+} := T_n^{e,f} - (\ell_0^e + \ell_0^f).$$

Recall the time-line construction of Section 2.2 using $m$ independent "time-lines" (one corresponding to each edge).

Now, fix arbitrary edges $e$ and $g$. By ignoring all of the time-lines except the ones corresponding to edges $e$ and $g$, one obtains the construction of the reinforced random walk under the law $\bar{P}^{1,1}$. Call this process the *restriction to edges $e$ and $g$*. Define

$$T^{e,g,restr} := \text{time of attraction for the restriction to } e \text{ and } g$$

and

$$T^{e,g,restr,+} := T^{e,g,restr} - (\ell_0^e + \ell_0^g).$$

In particular, $T^{e,g,restr,+}$ under $P^{\mathcal{G}}$ has the law of $T^+$ under $\bar{P}^{1,1}$. Next, observe that in the case where $g = f$ is the attracting edge, we have

(44)        $T^{e,f,restr} \equiv T_{Y^{e,f}}^{e,f} + 1$   and   $T^{e,f,restr,+} \equiv T_{Y^{e,f}}^{e,f,+} + 1,$

where the extra 1 on the right-hand side accounts for the traversal of edge $f$ at the attraction time $T^{e,f,restr}$. In addition, note that for each $e \neq f$, the number $2K_e$ of steps traversing $e$ before time $T$ equals the number of steps traversing $e$ before time $T_{Y^{e,f}}^{e,f}$. By similar reasoning, the number of steps traversing $f$ strictly before time $T$ equals the number of steps traversing $f$ before time $T_{Y^{g,f}}^{g,f}$, for at least one $g \neq f$ (to be precise, $g$ is the edge traversed at time $T - 1$). The last two claims imply that

(45)        $$T^+ \leq \sum_{e:e \neq f} T_{Y^{e,f}}^{e,f,+} + 1,$$

where, again, the extra 1 accounts for the traversal of $f$ at time $T$. By (44),

$$\sum_{e:e \neq f} T_{Y^{e,f}}^{e,f,+} + 1 = \sum_{e:e \neq f} (T^{e,f,restr,+} - 1) + 1 \leq \frac{1}{2} \sum_{e,g:e \neq g} T^{e,g,restr,+},$$

so (45) implies that

(46)        $$T^+ \leq \frac{1}{2} \sum_{e,g:e \neq g} T^{e,g,restr,+},$$

in particular, yielding (43). As noted already, the $\binom{m}{2}$ different random variables $T^{e,g,restr,+}$ are (identically) distributed under the law $P^{\mathcal{G}}$ as $T^+$ is under the law $\bar{P}^{1,1}$. The statement of the lemma is now a standard consequence of (43).   □

Now, consider a finite tree $\mathcal{G}$. Let $m(\mathcal{G})$ be the total number of pairs of edges in $\mathcal{G}$ that meet at a vertex. For example, the star with $m$ fingers has $m(\mathcal{G}) = \binom{m}{2}$. Using the same reasoning (for each $v$, consider separately the star created by restricting the tree to $v$ and all $u$, $u \sim v$) as in Lemma 14, one quickly obtains the following.

LEMMA 15.   *Assume that* $\ell_0^e = 1$, $e \in E(\mathcal{G})$. *Then,*

$$P^{\mathcal{G}}(T^+ > k) \leq m(\mathcal{G})(\bar{P}^{1,1}(T > k/m(\mathcal{G})) + \bar{P}^{1,2}(T > k/m(\mathcal{G}))).$$

Here, $\bar{P}^{1,2}(T > k/m(\mathcal{G}))$ appears due to parity considerations. Namely, for any two edges $e$, $f$ that meet at a vertex $v$, say, at the first time the walk visits $v$, the configuration of weights is either 1, 1 or 1, 2 or 2, 1.

We will soon show analogous results for the walk on a general finite graph. Before this, we quickly turn to the case where $\mathcal{G}$ is an infinite tree of bounded degree. Recall that $\#\mathcal{G}_1$ denotes the total number of vertices ever visited by the edge-reinforced random walk on $\mathcal{G}$. Here, again, we assume that $\ell_0^e = 1$, $e \in E(\mathcal{G})$. The next lemma can be proven in an analogous (but simpler) way to Lemma 25; we leave its verification to an interested reader.

LEMMA 16. *There exists $p > 0$, depending only on the weight $w$ and the degree $D(\mathcal{G})$ of $\mathcal{G}$, such that $\#\mathcal{G}_1/2$ is stochastically bounded by $Z$, where $Z$ is a geometric random variable with success probability $p$.*

COROLLARY 17. *Let $\mathcal{G}$ be an infinite tree such that $D(\mathcal{G}) < \infty$. Then, for any $c > 1$, we have*

$$P^{\mathcal{G}}(T^+ > k) \le O\left(\frac{1}{k^{c\log(1/(1-p))/2}}\right) + \max_{\mathcal{G}_{k,c}} P^{\mathcal{G}_{k,c}}(T^+ > k),$$

*where the above maximum is taken over all trees $\mathcal{G}_{k,c}$ having fewer than $cD(\mathcal{G})\log k$ vertices and degree bounded by $D(\mathcal{G})$.*

PROOF. Due to the last lemma, with probability $(1-p)^{(c\log k)/2}$, the range $\mathcal{G}_1$ of the walk is a subtree of $\mathcal{G}$ containing initial position $I_{t_0}$ and $c\log k$ or more vertices. On the opposite event, denoted by $B_{c\log k}$, we have $\mathcal{G}_1 \subset \mathcal{G}_{c\log k}^*$, where $\mathcal{G}_{c\log k}^*$ is a (nonrandom) subtree of $\mathcal{G}$ generated by all vertices $v$ of $\mathcal{G}$ such that the graph distance of $v$ and $I_{t_0}$ is less than or equal to $c\log k$. Therefore, by Lemma 13, we can bound

$$P^{\mathcal{G}}(\{T^+ > k\} \cap B_{c\log k}) \le P^{\mathcal{G}}(\{T^+ > k\} \cap \{\mathcal{G}_1 \subset \mathcal{G}_{c\log k}^*\}) \le P^{\mathcal{G}_{c\log k}^*}(T^+ > k). \square$$

COROLLARY 18. *Let $\mathcal{G}$ be an infinite tree of bounded degree and let $w(\cdot)$ be as in the examples of Theorem 10(a)–(b). Then,*

$$P^{\mathcal{G}}(T > k) = P^{1,1}(T > k)O(\log^q k)$$

*for any $q > \rho - \rho'$.*

PROOF. For any tree $\mathcal{G}_{k,c}$ of bounded degree with fewer than $O(\log k)$ vertices, one also has $m(\mathcal{G}_{k,c}) = O(\log k)$. Use the previous corollary with $c\log(1/(1-p))/2 > \rho - \rho' - 1$. Finally, note that under the assumptions of Theorem 10(a) [resp., (b)], $m(\mathcal{G}_{k,c})P^{1,1}(T^+ > k/m(\mathcal{G}_{k,c}))$ is of order $\log k \cdot P^{1,1}(T^+ > k)(\log k)^{\rho-\rho'-1}$ [resp., $\log k \cdot P^{1,1}(T^+ > k)(\log k)^{\rho-\rho'-1}(\log\log k)^{\alpha}$]. $\square$

3.2. *Analysis on finite graphs.* Assume that (A0) holds. Let $\mathcal{G}$ be a finite graph and let $\bar{n} = |E(\mathcal{G})|$. Moreover, denote the edges of $\mathcal{G}$ by $E(\mathcal{G}) = \{e_1, e_2, \ldots, e_{\bar{n}}\}$. If $v$ is an arbitrary vertex of the graph, let $n_v = \text{degree}(v)$, and let $\mathcal{N}_v := \{e_1^v, e_2^v, \ldots, e_{n_v}^v\}$ be the set of edges incident to $v$. Recall that $X_k^e$ equals the initial weight $\ell_0^e$ incremented by the number of times that edge $e$ has been visited by time $k$.

As before, we will start the walk at time $\sum_{e \in E(\mathcal{G})} \ell_0^e$. Fix the initial position $I_{t_0}$ at some arbitrary vertex $v_0$. The following proposition then holds.

PROPOSITION 19.  *Let $k \geq \sum_{e \in E(\mathcal{G})} \ell_0^e$ and $v \in V(\mathcal{G})$, and denote by $A_{v,k}$ the event $\{I_k = v\}$. Then, for any $\ell^e$, $e \in E(\mathcal{G})$ such that $\ell^e \geq \ell_0^e$, $e \in E(\mathcal{G})$ and $\sum_{e \in E(\mathcal{G})} \ell^e = k$, we have*

$$(47) \quad P^{\mathcal{G}}(X_k^e = \ell^e, e \in E(\mathcal{G}), A_{v,k}) \leq \frac{\prod_{e \in E(\mathcal{G})} w(\ell_0^e)}{\min_{e \in \mathcal{N}_{v_0}} w(\ell_0^e)} \cdot \frac{\sum_{e \in \mathcal{N}_v} w(\ell^e)}{\prod_{e \in E(\mathcal{G})} w(\ell^e)}.$$

REMARK.  Inequality (47) holds trivially when the conditions of the propositions do not hold, since the left-hand side then equals 0.

PROOF OF PROPOSITION 19.  As in the two-edge case, we will use induction on $\sum_e \ell^e = k$ to prove the above inequality. The base of induction at the initial time $\sum_e \ell_0^e$ clearly holds since, when the left-hand side is 0, the right-hand side is positive and when the left-hand side is 1, the right-hand side is greater than 1.

Now, take $k > \sum_e \ell_0^e$ and consider the event on the left-hand side. For each $i = 1, 2, \ldots, n_v$, let $v_i \in V(\mathcal{G})$ be the neighbor of $v$ such that $e_i^v = \{v, v_i\}$. In order for the event $\{X_k^e = \ell^e, e \in E(\mathcal{G}), A_{v,k}\}$ to occur, we must have $I_{k-1} = v_i$ for some $v_i \sim v$ such that $\ell^{\{v,v_i\}} = \ell^{e_i^v} > \ell_0^{e_i^v}$ and, furthermore, we must have $\{I_{k-1}, I_k\} = e_i^v$. Therefore,

$$P^{\mathcal{G}}(X_k^e = \ell^e, e \in E(\mathcal{G}), A_{v,k})$$
$$= \sum_{i=1, \ell_i^{e_i^v} > \ell_0^{e_i^v}}^{n_v} P^{\mathcal{G}}(X_{k-1}^e = \ell^e, \forall e \neq e_i^v, X_{k-1}^{e_i^v} = \ell^{e_i^v} - 1, A_{v_i,k-1})$$
$$\times \frac{w(\ell^{e_i^v} - 1)}{w(\ell^{e_i^v} - 1) + \sum_{e \neq e_i^v \in \mathcal{N}_{v_i}} w(\ell^e)}.$$

Similarly to Propositions 2 and 3, the proof follows immediately by induction.  □

From now on, denote by

$$\bar{w}_0(\bar{n}) := \frac{\prod_{e \in E(\mathcal{G})} w(\ell_0^e)}{\min_{e \in \mathcal{N}_{v_0}} w(\ell_0^e)}$$

the constant (ensuring appropriate scale-invariant behavior with respect to $w$) from the above proposition.

Define $S_1(k) := 1/w(k)$, $k \geq 1$, and for each $n \geq 2$ and $k \geq n$, define

$$(48) \qquad S_n(k) := \sum_{\ell_1 + \ell_2 + \cdots + \ell_n = k} \frac{1}{w(\ell_1)w(\ell_2) \cdots w(\ell_n)},$$

where the indices $\ell_i$, $i = 1, \ldots, n$, in the above summation are all greater than or equal to 1. If $k < n$, simply set $S_n(k) := 0$. Then, note that for $k \geq n \geq 2$,

$$
\begin{aligned}
(49) \qquad S_n(k) = &\frac{1}{w(1)} S_{n-1}(k-1) + \frac{1}{w(2)} S_{n-1}(k-2) + \cdots \\
&+ \frac{1}{w(k-n+1)} S_{n-1}(n-1).
\end{aligned}
$$

Subsequently, we will make use of the following assumption on $w(k)$:

$$(A2) \qquad \sum_{i=1}^{k-1} \frac{1}{w(i)w(k-i)} \leq \frac{C_w}{w(1)} \cdot \frac{1}{w(k)}, \qquad k \geq 1,$$

where $C_w < \infty$ depends on $w(\cdot)$ up to scaling.

REMARK. The examples of Theorem 10(a)–(c) all satisfy (A2).

The next lemma will be useful in deriving Corollary 22 below.

LEMMA 20. *If* (A2) *holds, then, for all* $k \geq n \geq 2$,

$$S_n(k) \leq \frac{(C_w)^n}{(w(1))^n w(k)}.$$

PROOF. We prove the statement inductively. The case $n = 2$ is a direct consequence of assumption (A2). Suppose that for some $n > 2$ and for all $k \geq n$, we have $S_n(k) \leq \frac{(C_w)^n}{(w(1))^n w(k)}$. Then, assumption (A2) and identity (49) imply, together with the inductive hypothesis, that for each $k \geq n + 1$,

$$S_{n+1}(k) \leq \frac{(C_w)^n}{(w(1))^n} \sum_{j=1}^{k-n} \frac{1}{w(j)w(k-j)} \leq \frac{(C_w)^{n+1}}{(w(1))^{n+1}w(k)}. \qquad \square$$

The next result is in the spirit of Lemma 6. It applies in the following setting: fix three different vertices $\omega, v$ and $u$ such that $\omega \sim v$ and $v \sim u$. Recall the notation from the beginning of this section. Furthermore, we assume, without loss of generality, that

$$(50) \qquad e_1^v = e_1^u = \{u, v\}, \qquad e_1^\omega = e_2^v = \{\omega, v\}.$$

Assume that $n_\omega = q$, $n_v = p$ and $n_u = m$ (recall that these are the degrees of the corresponding vertices). We introduce the following notation, to be used in the next theorem:

$$P^{\mathcal{G}}(\omega, u, v; k) := P^{\mathcal{G}}(I_k = \omega, I_{k+2i+1} = v, I_{k+2i+2} = u, i \geq 0).$$

THEOREM 21. *In the setting of Proposition* 19, *we have*

$$P^{\mathcal{G}}(\omega, u, v; k)$$

$$\leq \bar{w}_0(\bar{n}) \cdot \sum_{\ell^e : \sum_e \ell^e = k} \frac{w(\ell^{e_2^v})}{\prod_e w(\ell^e)}$$

(51)
$$\times \prod_{i=0}^{\infty} w(\ell^{e_1^v} + 2i) / (w(\ell^{e_1^v} + 2i) + w(\ell^{e_2^v} + 1)$$

$$+ w(\ell^{e_3^v}) + \cdots + w(\ell^{e_p^v}))$$

$$\times \prod_{i=0}^{\infty} \frac{w(\ell^{e_1^v} + 2i + 1)}{w(\ell^{e_1^v} + 2i + 1) + w(\ell^{e_2^v}) + \cdots + w(\ell^{e_m^u})}.$$

PROOF. This is a direct consequence of Proposition 19 and repeated conditioning. Namely, given a particular configuration of weights $\ell^e, e \in E(\mathcal{G})$, (47) estimates the probability for the walk to realize this configuration at time $k$ and to end up at vertex $\omega$ at time $k$, the probability of the next step is

$$w(\ell^{e_1^\omega}) / (w(\ell^{e_1^\omega}) + \cdots + w(\ell^{e_q^\omega})) = w(\ell^{e_2^v}) / (w(\ell^{e_1^\omega}) + \cdots + w(\ell^{e_q^\omega}))$$

and that of the infinitely many steps, each traversing $\{v, u\}$, is given by the two infinite products in the statement. Note that we have made use of the notation (50). □

There are various ways to simplify (and lose precision in doing so) the above bound. We chose a particularly simple one for the purposes of illustration since we could not find a good enough simplification that would "eliminate" the exponential term in the size $\bar{n}$ of the graph in Corollary 23 below. From now on, assume that both (A1) and (A2) hold.

Note that we can bound the sum (51) by

(52)
$$\sum_{\ell^e : \sum_e \ell^e = k} \frac{w(\ell^{e_2^v})}{\prod_e w(\ell^e)} \cdot \prod_{i=0}^{\infty} \frac{w(\ell^{e_1^v} + 2i)}{w(\ell^{e_1^v} + 2i) + w(\ell^{e_2^v} + 1)}.$$

Next, rearranging (52) according to the value $s = \ell^{e_1^v} + \ell^{e_2^v}$ yields [recall definitions (48) and (15)]

$$
\begin{aligned}
P^{\mathcal{G}}(\omega, v, u; k) & \\
& \leq \bar{w}_0(\bar{n}) \sum_{s=\ell_0^{e_1^v}+\ell_0^{e_2^v}}^{k-\sum_{e\neq e_1^v,e_2^v}\ell_0^e} S_{\bar{n}-2}(k-s) \sum_{j=1}^{s-1} \frac{1}{w(s-j)} \overline{W}_{s+1}(j+1) \\
& \leq \bar{w}_0(\bar{n}) \sum_{s=2}^{k-1} S_{\bar{n}-2}(k-s) \sum_{j=1}^{s-1} \frac{1}{w(s-j)} \overline{W}_{s+1}(j+1) \\
& \leq \bar{w}_0(\bar{n}) \sum_{s=2}^{k-1} S_{\bar{n}-2}(k-s) \sum_{j=1}^{s-1} \frac{1}{w(s-j)} \overline{W}_{k+1}(j+1),
\end{aligned}
\tag{53}
$$

where, for the very last inequality, we used (A1), which implies $\overline{W}_{s+1}(j+1) \leq \overline{W}_{k+1}(j+1)$, as in Lemma 7(b). Interchanging the order of summation, applying Lemma 20 and (A2) now gives

$$
P^{\mathcal{G}}(\omega, v, u; k) \leq \frac{\bar{w}_0(\bar{n})(C_w)^{\bar{n}-2}}{(w(1))^{\bar{n}-2}} \sum_{j=1}^{k-1} \sum_{s=j+1}^{k-1} \frac{1}{w(k-s)} \frac{1}{w(s-j)} \overline{W}_{k+1}(j+1)
$$

and, in turn,

$$
P^{\mathcal{G}}(\omega, v, u; k) \leq \frac{\bar{w}_0(\bar{n})(C_w)^{\bar{n}-1}}{(w(1))^{\bar{n}-1}} \sum_{j=1}^{k-1} \frac{1}{w(k-1-(j-1))} \overline{W}_{k+1}(j+1),
$$

which, comparing with the expression for $\bar{P}^{a,b}$ in Lemma 6 and accounting for various possibilities of parity, finally implies that

$$
P^{\mathcal{G}}(\omega, v, u; k) \leq \bar{w}_1 \frac{\bar{w}_0(\bar{n})(C_w)^{\bar{n}-1}}{(w(1))^{\bar{n}-1}} \sum_{a,b\in\{1,2\}} \bar{P}^{a,b}(T = k+1),
$$

where $\bar{w}_1 \in (0, \infty)$ accounts for the "$\asymp$" equivalence of Lemma 6.

REMARK. It will be convenient for the comparison arguments in the next corollary to refer to $\bar{P}^{a,b}$, even when $a$ and $b$ are of the same parity. In this case, the reader has an option of either noting that the parity does not influence the arguments for Theorem 4, Lemma 6 and Theorem 10, or noting that if $a - b$ is an even number, then the law $\bar{P}^{a,b}$ with reinforcement weight $w$ corresponds to the law $P^{a,b}$ with reinforcement weight $\bar{w}$, where $\bar{w}(a+j) = w(a+2j)$, $j \geq 0$.

COROLLARY 22. *Assuming* (A0)–(A2), *we have*

$$P^{\mathcal{G}}(T=k+1) \leq \bar{w}_1 \frac{2\bar{n}D(\mathcal{G})\bar{w}_0(\bar{n})(C_w)^{\bar{n}-1}}{(w(1))^{\bar{n}-1}} \sum_{a,b\in\{1,2\}} \bar{P}^{a,b}(T=k+1).$$

PROOF. Sum over all possible choices of vertex $v$, edge $\{u,v\}$ and neighbor $\omega$ of $v$, and note that there are at most $2\bar{n}D(\mathcal{G})$ terms of type $P^{\mathcal{G}}(\omega,v,u;k)$ contributing. □

COROLLARY 23. *Assuming* (A0)–(A2) *and* $\ell^e = 1$, *for all* $e \in E(\mathcal{G})$, *we have*

$$P^{\mathcal{G}}(T=k+1) \leq 2\bar{w}_1\bar{n}D(\mathcal{G})(C_w)^{\bar{n}-1} \sum_{a,b\in\{1,2\}} \bar{P}^{a,b}(T=k+1).$$

As noted earlier, the examples (a)–(c) from Theorem 10 satisfy (A0)–(A2). In particular, Corollaries 11 and 22 now imply the following.

COROLLARY 24. *Let* $\mathcal{G}$ *be a finite graph. Suppose that* $w(k)$ *is as in* (a) *or* (b) *in Theorem* 10. *Then,* $E^{\mathcal{G}}(T)$ *is infinite if* $\rho \leq 1 + \frac{1+\sqrt{5}}{2}$ *and finite if* $\rho > 1 + \frac{1+\sqrt{5}}{2}$.

The examples (d)–(e) of Theorem 10 do not satisfy (A2). Here, one could use (53) with separately derived bounds on $S_n(k)$ to obtain bounds on $P^{\mathcal{G}}(T=k+1)$, as in Corollary 22. In particular, if $w(k) = e^{k^{\varepsilon}}$, $\varepsilon \in (0,1]$, then $S_n(k) \leq \binom{k-1}{n-1}/w(k)$ and the above reasoning, together with Theorem 10(d), gives

$$(54) \qquad P^{\mathcal{G}}(T=k+1) = O(k^{\bar{n}}) \frac{k^{1-\varepsilon}}{e^{(k/2)^{\varepsilon}}},$$

which is $\bar{P}^{1,1}(T=k+1)$ up to a polynomial correction.

3.3. *Extensions to bounded degree graphs.* Let $\mathcal{G}$ be an infinite graph of bounded degree and, as usual, let assumption (A0) hold. We wish to estimate

$$P(\#\mathcal{G}_1 > k),$$

where we recall that $\#\mathcal{G}_1$ denotes the number of vertices in the range of the walk. Since $D(\mathcal{G}) < \infty$, note that the above estimate will imply an estimate on $P(|\mathcal{G}_1| > k)$.

LEMMA 25. *The random variable* $\#\mathcal{G}_1$ *is stochastically bounded by* $2 \cdot Z$, *where* $Z$ *has geometric distribution with success probability* $p \in (0,1)$, *where* $p$ *depends only on* $w(\cdot)$, $D(\mathcal{G})$ *and the initial configuration of weights* $\ell_0^e, e \in E(\mathcal{G})$.

As a consequence, we obtain that whenever $\mathcal{G}$ has bounded degree, both $\#\mathcal{G}_1$ and $|\mathcal{G}_1|$ have exponential tails.

PROOF OF LEMMA 25. We will construct a coupling of $\mathcal{G}_1$ and $\mathcal{G}_1^*$ such that $\mathcal{G}_1 \subset \mathcal{G}_1^*$, almost surely, and such that the claim of the lemma holds for $\mathcal{G}_1^*$. Denote by $T_v$ the time of the first visit to the vertex $v$ of $\mathcal{G}$, where $T_v$ is infinite if the walk never visits $v$. If $\{T_v = n\}$, then either:

(i) at least one neighbor $v'$ of $v$ has not been visited by the walk before time $n$; or

(ii) all the neighbors of $v$ were visited by the walk up to time $n$.

First, suppose that case (i) happens. Then, "add" to $\mathcal{G}_1^*$ both vertices $v$ and $v'$, as well as the just-traversed edge leading to $v$ and the edge $\{v, v'\}$. Due to the assumptions, with probability $p$, uniformly bounded away from 0 and depending only on $D(\mathcal{G})$, $w(\cdot)$ and $\ell_0^e, e \in E(\mathcal{G})$ (in fact only $\ell_0^f$ on edges $f$ incident to $v'$), the walk keeps traversing solely the edge $\{v, v'\}$ after time $n$. In symbols, on the event of case (i),

$$P(\{I_k, I_{k+1}\} = \{v, v'\}, k \geq n | \mathcal{F}_n) > p.$$

If the above event $\{I_k, I_{k+1}\} = \{v, v'\}, k \geq n$, does not occur, then the walk will keep exploring the graph elsewhere. Either it will get attracted to an edge before encountering another new vertex or it will encounter another new vertex prior to getting attracted.

If the case (ii) happens, note that $\mathcal{G}_1^*$ already contains vertex $v$. Namely, let $u$ be the neighbor of $v$ such that

$$T_u = \max_{v' \sim v} T_{v'} < T_v.$$

Then, $v$ must have been added to $\mathcal{G}_1^*$ as part of the case (i) procedure, before or at time $T_u$.

Therefore, $\mathcal{G}_1 \subset \mathcal{G}_1^*$, almost surely, by induction. Moreover, from the construction, it is clear that $\#\mathcal{G}_1$ is stochastically bounded by $2Z$. □

As a conclusion, we offer the following weak universality-type result.

COROLLARY 26. *Assume that $\ell^e = 1$ for all $e \in E(\mathcal{G})$. If $w(k) = k^\rho$, $\rho > 1$, then there exists $p > 0$ such that $P^{\mathcal{G}}(T > k) \leq \frac{1}{k^p}$.*

PROOF. We use the idea of Corollary 17, together with the previous lemma and the bound of Corollary 23. Namely, we split the event $\{T > k\}$ according to whether or not the walk reaches distance $d_k$ from $I_0$. Choose $d_k = c \log k$, where $c$ is such that $(C_w)^{c \log k} \ll k^{\rho - \rho' - 1}$. □

REMARK.    If $w(k) = e^k$, then the walk gets attracted at any particular step with probability bounded away from 0, so there exists $c > 0$ such that $P^{\mathcal{G}}(T > k) = O(\frac{1}{e^{ck}})$. Somewhat disappointingly, the bound of type (54) is too weak to provide an alternative derivation (analogous to the proof of the last corollary) of the above bound. Indeed, the question of finding the exact (up to a multiplicative constant) behavior of the tail distribution of $T$ on general bounded degree graphs, even in the case of the examples in Theorem 10, remains open.

**Acknowledgments.**    Research conducted while holding a postdoctoral fellow position at the University of British Columbia, Vancouver. The first author is very grateful to David Brydges for his advice and support during her time as a postdoctoral fellow at the University of British Columbia. We wish to thank the anonymous referee for many helpful comments and suggestions.

TU BERLIN
FAKULTÄT II
INSTITUT FÜR MATHEMATIK
STRASSE DES 17 JUNI 136
D-10623 BERLIN
GERMANY
E-MAIL: cotar@math.tu-berlin.de

UNIVERSITÉ DE PROVENCE
LATP UMR 6632
CENTRE DE MATHÉMATIQUES ET INFORMATIQUE
39 RUE DE F. JOLIOT-CURIE
13453 MARSEILLE, CEDEX 13
FRANCE
E-MAIL: vlada@cmi.univ-mrs.fr